\documentclass[12pt]{amsart}
 \usepackage{amsfonts,amssymb,eucal}
\usepackage{amsthm} 
 \usepackage{amsmath} 
\usepackage{amscd}
 \usepackage{latexsym} 
\input{cyracc.def} 
\font\tencyr=wncyr10
\def\cyr{\tencyr\cyracc}
\numberwithin{equation}{section}
\newcommand{\Ad}{\ensuremath{{\mbox{\rm{Ad}}}}}
\newcommand{\ad}{{\mbox{\rm ad}}}
\newcommand{\e}{{\mbox{\rm e}}}

\newcommand{\mb}[1]{{\mbox{\boldmath{$#1$}}}}
\newcommand{\mc}[1]{{\mathcal{#1}}}
\newcommand{\got}[1]{{\mathfrak{#1}}}
\newcommand{\db}[1]{{\mathbb{#1}}}
\newcommand{\gata}{\blacksquare}
\newcommand{\pa}{\partial}
\newcommand{\R}{\ensuremath{\mathbb{R}}}
\newcommand{\C}{\ensuremath{\mathbb{C}}}
\newcommand{\N}{\ensuremath{\mathbb{N}}}

\newcommand{\Hi}{\ensuremath{\mathcal{H}}}

\newcommand{\Hinf}{\ensuremath{\mathcal{H}^{\infty}}}
\newcommand{\g}{\ensuremath{\got{g}}}
\newcommand{\gc}{\ensuremath{\got{g}_{\C}}}
\newcommand{\Ugc}{\ensuremath{\mathcal{U}({\gc}})}
\renewcommand{\P}{\ensuremath{\mathbb{P}}}
\newcommand{\Ph}{\ensuremath{\P (\Hi )}}
\newcommand{\Phinf}{\ensuremath{\P (\Hinf )}}
\newtheorem{Remark}{Remark}

\newcommand{\DM}{\ensuremath{{\got{D}}_M }}
\newcommand{\AM}{\ensuremath{{\got{A}}_M }}
\newcommand{\D}{\ensuremath{{\got{D}}}}
\newcommand{\AAA}{\ensuremath{{\db{A}}_M}}
\newcommand{\am}{\ensuremath{{\bf{A}}_M}}%

\newcommand{\FSB}{symmetric Fock space }

\newcommand{\fl}{\ensuremath{{\mathcal{F}}_{\Hi}}}
\newcommand{\NC}{\ensuremath{\mathcal{G}}}

\newtheorem{Proposition}{Proposition}
\newtheorem{lemma}{Lemma}
\newtheorem{corollary}{Corollary}
\newtheorem{Comment}{Comment}
\theoremstyle{definition}

 \newcommand{\ka}{\ensuremath{{\kappa}}}
\newcommand{\SL}{\text{SL}_2(\R )}
\begin{document}

\title{A holomorphic representation of the Jacobi algebra}
\author{Stefan  Berceanu}
\address[Stefan  Berceanu]{National
 Institute for Physics and Nuclear Engineering\\
         Department of Theoretical Physics\\
         PO BOX MG-6, Bucharest-Magurele, Romania}
\email{Berceanu@theory.nipne.ro}

\begin{abstract} A representation of the Jacobi algebra
$\mathfrak{h}_1\rtimes \mathfrak{su}(1,1)$ by first order differential
operators with polynomial coefficients on the manifold $\mathbb{C}\times
\mathcal{D}_1$
 is presented. The Hilbert space of holomorphic functions
 on which the holomorphic
first order differential operators with polynomials coefficients
 act is constructed.
\end{abstract}
\subjclass{81R30,32WXX,12E10,33C47,32Q15,81V80}
\keywords{Coherent states,
representations of coherent state Lie algebras, Jacobi group,
first order  holomorphic differential
operators with polynomial coefficients}
\maketitle
\noindent
\tableofcontents
\section{Introduction}
In this paper we deal with realizations of finite-dimensional Lie
algebras by first-order differential operators on homogeneous
spaces. Our method, firstly developed in   \cite{morse},
 permits to get the holomorphic differential action of the generators
of a continuous unitary representation $\pi$
   of a Lie  group $G$ with the Lie algebra $\got{g}$  on a homogeneous
space $M=G/H$. We consider 
homogeneous manifolds realized as K\"ahler coherent
state (CS)-orbits  obtained by the action of the 
representation
$\pi$ on a fixed cyclic vector $e_0$ belonging to the complex separable
Hilbert
space $\Hi$   of the representation \cite{perG}.
 We have applied our  method to compact (non-compact)
hermitian symmetric spaces in \cite{sbcag} (respectively, \cite{sbl})
and we have produced simple formulas which show that the differential
action of the generators of a hermitian group $G$ on holomorphic
functions defined on  the hermitian
symmetric spaces $G/H$ can be written down as a sum of two terms, one
a polynomial $P$, and the second one a sum of partial derivatives times
some polynomials $Q$-s, the degree of polynomials being less than
3. This  is a generalization of the well  known realization \cite{lie} of the
generators $J_{0,+,-}$ 
of the group $G=\text{SU}(2)$ (and similarly
for its non-compact
dual $G=\text{SU}(1,1)$) on the homogeneous manifold  $G/\text{U}(1)$ by the
differential operators $\db{J}_+=-\frac{\pa}{\pa z}$,
$\db{J}_-=-2jz+z^2\frac{\pa}{\pa z}$, $\db{J}_0=j-z\frac{\pa}{\pa
z }$,  where the generators verify the commutation relations
$[J_0,J_{\pm}]=\pm J_{\pm}$, $[J_-,J_+]=-2J_0$ and $J_0e_0=-je_0$.

In \cite{sbctim}, \cite{sin} we have generalized the results
 of \cite{sbcag}, \cite{sbl} to
K\"ahler CS-orbits of semisimple Lie groups. The
differential action of the generators of the
 groups is  of the same type  as in the
case of hermitian symmetric orbits,  i.e. first order differential
operators with holomorphic polynomial coefficients, but the maximal  degree of
the polynomials is  grater than 2. We have
presented explicit formulas involving the Bernoulli numbers and the
structure constants for semisimple Lie groups \cite{sbctim},
\cite{sin}.
 The simplest example in which the
maximum degree of the polynomials multiplying the derivative is already 3 was
worked out  in detail in \cite{sbctim}, \cite{sin}, where we have
   constructed CS  on 
the non-symmetric space $M:= \text{SU}(3)/\text{S}(\text{U}(1)\times \text{U}(1)\times \text{U}(1))$.

Let us now 
recall  the standard  Segal-Bargmann-Fock \cite{bar}
 realization
$a\mapsto \frac{\partial}{\partial z};~~ a^+\mapsto z$ of the
canonical commutation relations (CCR) $[a,a^+]=1$  on the
symmetric Fock space $\fl :=\Gamma^{\mbox{\rm{hol}}}(\C,\frac{i}{2
\pi}\exp(-|z|^2)dz\wedge d\bar{z})$  attached to the Hilbert space
$\Hi :=L^2(\R, dx)$.
The Segal-Bargmann-Fock realization can be considered as
a representation by differential operators of the real 3-dimensional
 Heisenberg algebra
  $\got{h}_1 \equiv \g_{HW}=
< is1 + z{\mb{a}}^+ -\bar{z}{\mb{a}} >_{s\in\R;z\in\C}$ of
 the Heisenberg-Weyl group (HW) $H_1$, where $H_n$ denotes the $(2n+1)$-dimensional
HW group. We can look at this construction from
 group-theoretic point of view,
considering  the complex number $z$
 as local coordinate on the homogeneous manifold
 $M:=H_1/\R\cong \C$.  Glauber \cite{gl}
 has attached
 field coherent states  to the points of the manifold $M$.

In the present paper we are interested in  representations of Lie
algebras which are semi-direct sum  of
Heisenberg   algebras and  semisimple Lie algebras 
 by first order differential operators with holomorphic
polynomials coefficients. The most appropriate framework for such
an approach   is furnished by the so called CS-groups, i.e. groups which
admit an orbit which is a complex  submanifold of a projective Hilbert
space \cite{lis},\cite{neeb}. Indeed, such groups contain all compact
groups,
all simple hermitian groups, certain solvable groups and also some
mixed groups as the semi-direct product of the HW group and the
symplectic group \cite{neeb}. In reference \cite{last} we have
advanced the hypothesis that {\it the generators
of  CS-groups admit representations by 
first order differential operators with holomorphic
polynomials coefficients on CS-manifolds}. Here we just present
explicit formulas for the
simplest example of such a  representation of the 
Lie  algebra
semi-direct sum  of the three-dimensional Heisenberg   algebra
 $\got{h}_1$ and the
algebra of the group $\text{SU}(1,1)$ acting on it in the canonical
fashion,
 $\got{g}^J_1:=\got{h}_1\rtimes \got{su}(1,1)$,
called {\it Jacobi algebra} (cf. \cite{ez} or p. 178  in
\cite{neeb}). The case of the  
Jacobi algebra $\got{g}^J_n=\got{h}_n\rtimes \got{sp}(n,\R )$ is
treated  separately \cite{sbj}. Let us remained also that 
the Jacobi algebra $\got{g}^J_n$, also denoted $ \got{st}( n,\R )$
 by Kirillov
 in \S 18.4 of \cite{kir} or $\got{tsp}(2n+2,\R )$ in \cite{kir2},
 is isomorphic with the subalgebra of Weyl algebra  $A_n$ (see also
\cite{dix}) of polynomials of degree maximum 2 in the variables
$p_1,\dots,p_n, q_1,\dots,q_n$ with the Poisson bracket, while  the
 Heisenberg algebra $\got{h}_n$ is the nilpotent ideal isomorphic with
polynomials 
of degree  $\le 1$ and  the real symplectic algebra $\got{sp}(n,\R
)$ is isomorphic to the subspace of symmetric homogeneous polynomials
of degree $2$.  In this paper we study the 6-dimensional 
Jacobi algebra $\got{g}^J_1$ and we denote it just $\got{g}^J$, when
there is no possibility of confusion with $\got{g}^J_n$. 

The representations of the Jacobi group were investigated also by the
orbit method \cite{kir,kir1}, starting from a matrix representation
(see p. 182 in \cite{kir1}) of the
Jacobi algebra $\got{g}^J$ in
\cite{bs,berndt}. Our method is inspired from the squeezed states of
Quantum Optics, see e.g. the reviews 
 \cite{zh,siv,dod,dr}. It is well known that for the
harmonic oscillator CSs the uncertainties in momentum and position
are equal with $1/\sqrt{2}$ (in units of $\hbar$). ``The squeezed
states'' \cite{ken,stol,lu,bi,yu,ho,wa}
are the states for which the uncertainty in  position is less than
$1/\sqrt{2}$. The squeezed states are a particular class of ``minimum
uncertainty states'' (MUS) \cite{mo}, 
i.e. states which saturates the Heisenberg
uncertainty relation. 
 In the present paper we do not
insist on the applications of our paper to  the squeezed states, the
Gaussian states \cite{si,ali},  disentangling theorems,
i.e. analytic Backer-Campbell-Hausdorff relations defined from a
$4\times4$-matrix representation of the Jacobi algebra, 
  or nonlinear coherent states
\cite{siv}. Let us just mention that ``Gaussian pure states''
(``Gaussons'') \cite{si} are more general MUSs. In fact, as was shown in
\cite{ali}, these states are CSs based on the manifold
$\mc{X}^J_n:= \mc{H}_n\times\R^{2n}$, where $\mc{H}_n$ is the Siegel upper half
plane  $\mc{H}_n:=\{Z\in M_n(\C)| Z=U+iV, U,V\in M_n(\R), \Im (V)>0, U^t=U;
V^t=V\}$.
$M_n(R)$ denotes the $n\times n$ matrices with entries in $R$,
$R=\R~\text{or}~\C$ and $X^t$ denotes the transpose of the matrix
$X$. In  \cite{sbj}
 we have started the generalization of CSs attached to the
 Jacobi group $G^J_1=H_1\rtimes\text{SU}(1,1)$
 to the Jacobi group
$G^J_n=H_n\rtimes\text{Sp}(n,\R )$. The connection of  our
construction of coherent states based on $\mc{D}^J_n=\C^n\times\mc{D}_n$
\cite{sbj}  and
the Gaussons of \cite{si} is a subtle one and should be investigated
separately. $\mc{D}_n$ denotes the Siegel ball  $\mc{D}_n:= \{Z\in
M_n(\C)|Z=Z^t, ~ 1-Z\bar{Z}> 0 \}$. In \S \ref{rem19} we  
indicate the clue  of this connection
  in the present case,
$n=1$,  which is offered by the K\"ahler-Berndt's construction, shortly 
sketched  in the same section \S \ref{rem19}.
      The only physical applications are contained in \S
\ref{app}, where we use the expressions of the generators of the
Jacobi group $G^J_1$ to determine the quantum and classical evolution
on the manifold $\mc{D}^J_1$,   
generated by a linear Hamiltonian  in the generators of the group. 

We  emphasize 
 that some of the results obtained in this paper, as the
reproducing kernel or the group action on the base manifold, can be
obtained as particular cases of some of the formulas in Chapter III,
 Propositions
5.1-5.3 in \cite{satake} and
\S XII.4 in \cite{neeb}. We also stress that some of the formulas presented
here appear in the context of  automorphic Jacobi forms \cite{ez,bs} -this
denomination is inspired by the book \cite{ps}. The
Jacobi group can be associated (see  Chapter 5  entitled {\it ``K\"ahler's New
Poincar\'e group''} in the article {\it ``Survey of K\"ahler's mathematical
work and some comments''} of R. Berndt and O. Riemenschneider in
\cite{cal}) with the group $G^K$ investigated by K\"ahler
\cite{cal1,cal2,cal3} 
as a group of the
{\it Universal Theory of Everything}, including relativity, quantum
mechanics and even biology.  In the paper \cite{cal1} K\"ahler  has
determined  the structure of the 
real 10-dimensional Lie algebra $\got{g}^K$ of the ({\it Poincar\'e}
or {\it New Poincar\'e}) group  $G^K$ and has realized this algebra by
differential operators in four {\bf  real} variables.  
  However, our approach and the  proofs are independent and, we hope, more
accessible to people familiar with the coherent state approach in
Theoretical Physics and in 
Mathematical Physics. Moreover, as far as we know, some of the formulas
presented in this paper
are completely new, e.g.  (\ref{x3x})
expressing  the base of polynomials defined on the 
manifold $\mc{D}^J_1$ -
 the homogeneous space of the Jacobi group $G^J_1$, acting by
biholomorphic maps, or the resolution of unity (\ref{ofi})-(\ref{ofi3}).  

In
order to facilitate the understanding of all  subsequent sections, we
present in \S \ref{css}  the general setting concerning the CS-groups:
\S \ref{CSgr} briefly recalls the definition of CS-groups and \S \ref{s22}
defines the space of functions, called the symmetric Fock  space,  on
which the the differential operators  act (\S \ref{s33}).
 However, we shall not enter into a detailed analysis of the
root structure of CS-Lie algebras \cite{neeb},
 keeping the exposition as elementary
as possible. \S \ref{jac1} presents the Jacobi algebra $\got{g}^J$. Perelomov's
CS-vectors associated with the Jacobi group $G^J_1$ (cf.
denomination used in  \cite{bs} or 
at p. 701 in \cite{neeb}) are based on the complex homogeneous 
manifold $M:= \mc{D}^J_1$.  The differential
action of the generators of the Jacobi group is given  in Lemma
\ref{mixt} of \S \ref{diff}.  The operators $a$ and $a^+$ are
unbounded operators, but it is enough to work on the dense subspace of
smooth vectors of the Hilbert space of the hermitian representation
(cf. p. 40 in \cite{neeb}  and also \S \ref{s33} of our paper).
   In Lemma \ref{lema5} of \S \ref{repK}
we calculate the reproducing kernel 
$K:\mc{D}^J_1\times \bar{\mc{D}}^J_1\rightarrow \C$. 
 Some facts concerning the representations of the
HW group $H_1$ and $\text{SU}(1,1)$ are collected in \S \ref{unul}. Several relations
are obtained in \S \ref{hpb} as a consequence of the fact that  
the Heisenberg  
 algebra is an ideal of the
Jacobi algebra, and we find  how  to change the order of the representations
of the groups HW and $\text{SU}(1,1)$.  Some of the relations presented  in
\S \ref{hpb}
have appeared earlier  in connection with the squeezed states \cite{ken} in
Quantum Optics \cite{stol}.
The main result of 
\S \ref{actg} is
given in Proposition \ref{mm1}, which expresses the action of the Jacobi
group
 on  Perelomov's CS-vectors. Remark \ref{rem9} establishes  the connection
of our results in the context of coherent states with  those obtained in
the theory
of automorphic Jacobi  forms \cite{ez}.  In \S \ref{jcg} we 
 construct the symmetric Fock space attached to the reproducing kernel
$K$   from the symmetric Fock spaces associated with the groups 
 HW (cf. \S \ref{heil})
and $\text{SU}(1,1)$ (cf. \S \ref{su1}). The $G^J_1$-invariant 
K\"ahler two-form $\omega$, the Liouville form and the equations of
geodesics
 on the manifold $\mc{D}^J_1$ are
calculated  in \S \ref{two}.
   Proposition \ref{final} summarizes all the
information obtained in  \S \ref{sFs} concerning  the symmetric Fock
space $\mc{F}_K$
attached to the reproducing kernel $K$  for   the Jacobi group
$G^J_1$, while  Proposition \ref{finalf} gives the continuous unitary
holomorphic representation $\pi_K$ of  $G^J_1$ on $\mc{F}_K$.
Simple applications to
equations of motion on $\mc{D}^J_1$ determined by linear Hamiltonians in the
generators of the Jacobi group are presented in \S \ref{app}. The
equation of motion is a matrix Riccati equation  on the
manifold $\mc{D}^J_1$.  In order to compare our  K\"ahler two-form $\omega$  with
that given by E. K\"ahler (see \cite{cal}, which reproduces
\cite{cal1,cal2,cal3})
and R. Berndt \cite{bern,bs}, we  express 
in \S \ref{rem19} 
 our $\omega$  in coordinates on $\mc{D}^J_1$ in
appropriate
(called in \cite{bs} {\it EZ}) coordinates
in $\mc{X}^J_1$.
The K\"ahler-Berndt's two-form is in fact the K\"ahler two-form
attached to the manifold $\mc{X}^J_1$ on which are based the
Gaussons considered in \cite{si} in the case n=1. \S \ref{lastc}
contains some more 
remarks referring to the
connection between  the
formulas proved in the present article for the Jacobi algebra
and the  formalism used in \cite{neeb} for CS-groups. 
In order to be self-contained, two  formulas referring to the groups HW and
$\text{SU}(1,1)$ are proved in the Appendix.

\section{Coherent states: the general setting}\label{css}
\subsection{Coherent state groups}\label{CSgr}

Let us consider the triplet $(G, \pi , \Hi )$, where $\pi$ is
 a continuous, unitary
representation
 of the  Lie group $G$
 on the   separable  complex  Hilbert space \Hi .
Let us denote by $\Hinf$ the {\it smooth vectors}.

Let us
pick up $e_0\in \Hinf$ and let  the notation:
$e_{g,0}:=\pi(g).e_0, g\in G$.
We have an action $G\times \Hinf\rightarrow\Hinf$, $g.e_0 :=
e_{g,0}$. When there is no possibility of confusion, we write just
$e_{g}$ for $e_{g,0}$. 

Let us denote by
$[~]:\Hi^{\times}:=\Hi\setminus\{0\}\rightarrow\Ph=\Hi^{\times}/ \sim$
the projection with respect to the equivalence relation
 $[\lambda x]\sim [x],~ \lambda\in \C^{\times},~x\in\Hi^{\times}$. So,
$[.]:\Hi^{\times}\rightarrow \Ph , ~[v]=\C v$. The action
$ G\times \Hinf\rightarrow\Hinf$ extends to the action
$G\times \Phinf\rightarrow\Phinf ,~ g.[v]:=[g.v]$.

For $X\in\g$, where \g ~is the Lie algebra of the Lie group $G$,  let
us define the (unbounded) operator $d\pi(X)$ on \Hi~ by
$d\pi(X).v:=\left. {d}/{dt}\right|_{t=0} \pi(\exp tX).v$,
whenever the limit on the right hand side exists. 
   We obtain a
representation of the Lie algebra \g~ on \Hinf , {\it the derived
representation}, and we denote
${\mb{X}}.v:=d\pi(X).v$ for $X\in\g ,v\in \Hinf$. Extending $d\pi$ by complex
linearity, we get a representation of the universal enveloping algebra
of the complex Lie algebra \gc~ on
the complex vector space \Hinf ,
$ d\pi:\mc{S}:=\Ugc\rightarrow B_0(\Hinf ) $. 
 Here $B_0(\Hi^0)\subset \mc{L}(\Hi )$, where $\Hi^0:=\Hinf$ denotes
the  subset   of linear operators  $A:\Hi^0\rightarrow \Hi^0$
which have  a formal adjoint (cf. p. 29 in \cite{neeb}). 

Let us now denote by $H$  the isotropy group $H:=G_{[e_0]}:=
\{g\in G|g.e_0\in\C e_0\}$.
We shall consider (generalized) coherent
 states on complex  homogeneous manifolds $M\cong
G/H$ \cite{perG}, imposing the restriction that $M$ be 
{\it  a complex submanifold of \Phinf }. In such a case the orbit $M$ 
is called a {\it
CS-manifold} and the groups $G$ which generate such
orbits are called   {\it CS-groups} (cf.   Definition XV.2.1 at p. 650
and Theorem XV.1.1 at p. 646 in  \cite{neeb}),
while their Lie algebras are called {\it CS-Lie algebras}.

 {\em The coherent vector
 mapping} is defined locally, on a coordinate neighborhood $\mc{V}_0$,
 $\varphi : M\rightarrow \bar{\Hi}, ~ \varphi(z)=e_{\bar{z}}$
(cf. \cite{last}),
where $ \bar{\Hi}$ denotes the Hilbert space conjugate to $\Hi$.
The  vectors $e_{\bar{z}}\in\bar{\Hi}$ indexed by the points
 $z \in M $ are called  {\it
Perelomov's coherent state vectors}. The precise  definition depends on
the root  structure of the CS-Lie algebras and we do not go into
the details here (see  \cite{last}), but only in \S
\ref{lastc}
we just specify the root structure according to \cite{neeb}
 in the case of the Jacobi algebra.

We use for the scalar product the convention:
$(\lambda x,y)=\bar{\lambda}(x,y)$, $ x, y\in\Hi ,\lambda\in\C $.

\subsection{The symmetric Fock space}\label{s22}

The space of holomorphic functions (in fact, holomorphic sections of a
certain $G$-homogeneous line bundle over $M$ \cite{neeb}, \cite{last})
\fl~ is defined as the set of square integrable 
functions 
 with respect to  the scalar product
\begin{equation}\label{scf}
(f,g)_{\fl} =\int_{M}\bar{f}(z)g(z)d\nu_M(z,\bar{z}),
\end{equation}
\begin{equation}\label{scf1}
d{\nu}_{M}(z,\bar{z})=\frac{\Omega_M(z,\bar{z})}{(e_{\bar{z}},e_{\bar{z}})}.
\end{equation}
Here  $\Omega_M$ is the normalized  $G$-invariant volume form
\begin{equation}
\Omega_M:=(-1)^{\binom {n}{2}}\frac {1}{n!}\;
\underbrace{\omega\wedge\ldots\wedge\omega}_{\text{$n$ times}}\ ,
\end{equation}
and the $G$-invariant  K\"ahler two-form $\omega$ on the $2n$-dimensional
manifold  $M$ is given by
\begin{equation}\label{twoform}
\omega (z)
 = i\sum_{\alpha ,\beta}
\NC_{\alpha,\beta}dz_{\alpha} \wedge
d\bar{z}_{\beta},~
\NC_{\alpha,\beta}(z)=\frac{\pa^2}{\pa z_{\alpha}\pa\bar{z}_{\beta}}\log
(e_{\bar{z}},e_{\bar{z}}).
\end{equation}

It can be shown  that    (\ref{scf})
is nothing else that  {\em {Parseval   overcompletness
 identity}}  \cite{berezin}
 \begin{equation}\label{orthogk}
  (\psi_1,\psi_2)=\int_{M=G/H}(\psi_1,e_{\bar{z}})(e_{\bar{z}},\psi_2)
 d{\nu}_{M}(z,\bar{z}),~ (\psi_1,\psi_2\in \Hi ) .
 \end{equation}

It can be seen 
that {\it  relation  {\em (\ref{scf})} (or {\em  (\ref{orthogk})})
 on homogeneous
manifolds fits into
 Rawnsley's global realization \cite{raw} of Berezin's coherent states on
quantizable K\"ahler manifolds \cite{berezin}},
modulo   Rawnsley's ``epsilon'' function \cite{raw,CGR1},
 a constant for homogeneous quantization. 
If $(M,\omega )$ is a K\"ahler manifold and $(L,h, \nabla )$ is a (quantum)
holomorphic line bundle $L$ on $M$, where 
$h$ is the hermitian metric and $\nabla$
is the connection compatible with the metric and the complex structure,
then $h(z,\overline{z})=(e_{\bar{z}},e_{\bar{z}})^{-1}$ and the
K\"ahler potential $f$ is $f=-\log h(z)$ (see e.g. \cite{sbm}).

Let us  now introduce the map
\begin{equation}\label{aa}
\Phi :\Hi^{\star}\rightarrow \fl ,\Phi(\psi):=f_{\psi},
f_{\psi}(z)=\Phi(\psi )(z)=(\varphi (z),\psi)_{\Hi}=(e_{\bar{z}},\psi)_{\Hi},~
z\in{\mathcal{V}}_0,
\end{equation}
where we have identified the space $\overline{\Hi}$  complex conjugate
 to \Hi~  with the dual
space
$\Hi^{\star}$ of $\Hi$.

It can be defined a function 
   $K: M\times\overline{M}\rightarrow \C$, which on  ${\mathcal{V}}_0\times
 \overline{\mathcal{V}}_0$ reads
\begin{equation}\label{kernel}
K(z,\overline{w}):=K_w(z)=
(e_{\bar{z}},e_{\bar{w}})_{\Hi}.
\end{equation}

For CS-groups, {\it the function $K$ {\em(\ref{kernel})} is a positive
definite  reproducing
kernel; the \FSB~\fl~ (or $\mc{F}_K$)
 is the reproducing kernel Hilbert space of holomorphic functions on
$M$,
$\Hi_K\subset \C^M$, 
associated to the kernel $K$} (\ref{kernel}), and
 {\it the evaluation  map $\Phi$ defined in
 {\em  (\ref{aa})}
extends to an   isometric $G$-equivariant embedding
$\Hi^*\rightarrow\fl$}
 \cite{last}
\begin{equation}\label{anti}
(\psi_1,\psi_2)_{\Hi^{\star}}=(\Phi (\psi_1),\Phi
(\psi_2))_{\fl}=(f_{\psi_{1}},f_{\psi_{2}})_{\fl}=
\int_M\overline{f}_{\psi_1} (z)f_{\psi_2}(z)d\nu_M(z) .
\end{equation}
Sometimes the kernel $K$ is considered as a {\it Bergman section} \cite{PWJ} of a
certain bundle over $M\times\bar{M}$, firstly considered by Kobayashi
\cite{koby}, see Chapters V-VIII in \cite{neeb94} and Chapter XII in \cite{neeb}.

\subsection{Representation of CS-Lie algebras by differential
operators}\label{s33}

Let us consider again the triplet $(G, \pi, \Hi )$.
  {\it The derived representation}
 $d\pi$  is a {\it hermitian representation}
 of the semi-group
 $\mathcal{S}:=\Ugc$ on $\Hinf$ (cf. p. 30 in \cite{neeb}). The
unitarity and the continuity of  the representation $\pi$ 
imply that $i d\pi (X)|_{\Hinf}$ is essentially selfadjoint
(cf. p. 391 in \cite{neeb}).
Let us denote his image in $B_0(\Hinf )$ by $\am := d\pi(\mathcal{S})$.
If $\Phi : \Hi^*\rightarrow \fl $ is the isometry
(\ref{aa}), we are interested in the study of the image of
\am~  via $\Phi$  as subset in the
algebra of holomorphic, linear differential operators,
$ \Phi\am\Phi^{-1}:={\db{A}}_M\subset\got{D}_M$.

The   set
 $\got{D}_M$ (or simply \D ) {\it of holomorphic, finite order, linear
differential operators on} $M$ is a
 subalgebra of homomorphisms ${\mc Hom}_{\C}({\mc O}_M,{\mc O}_M)$
 generated
 by the set ${\mc O}_M$ of germs of holomorphic functions of $M$ and the
 vector fields.
 We consider also {\it  the subalgebra} \AM~ of ${\db{A}}_M$~
 {\it of differential operators with  holomorphic polynomial coefficients}.
Let $U:=\mc{V}_0\subset M$, endowed with the local coordinates
$(z_1,z_2,\cdots ,z_n)$. We set $\pa_i :=\frac{\pa}{\pa z_i}$ and
$\pa^{\alpha}:=\pa^{\alpha_1}_1
\pa^{\alpha_2}_2\cdots \pa^{\alpha_n}_n$, $\alpha :=(\alpha_1,
\alpha_2 ,\cdots ,\alpha_n)\in\N^n$. The sections of \DM~ on $U$ are
$A:f\mapsto \sum_{\alpha}a_{\alpha}\pa^{\alpha}f$,
$a_{\alpha}\in\Gamma (U, {\mc{O}})$,  $a_{\alpha}$-s being zero
except a finite number.

For $k\in\N$, let us denote by $\D_k$ the subset of differential
operators of degree $\le k$.  The filtration of
$\D$ induces a filtration on $\got{A}$. 

Summarizing, we have a correspondence between the following three objects:
\begin{equation}\label{correspond}
\g_{\C} \ni X \mapsto\mb{X}\in\am\mapsto\db{X}\in\AAA\subset \DM, {\mbox{\rm
 ~ {differential~operator~on}}}~ \fl .
\end{equation}

Moreover, it is easy to see \cite{last} that
 {\it if $\Phi$ is the isometry
{\em{(\ref{aa})}}, then
$\Phi d\pi(\g_{\C})\Phi^{-1}$ $\subseteq \D_1$ and we have}
\begin{equation}\label{car1}
\g_{\C}\ni X \mapsto{\db{X}}\in \D_1;~~
 {\db{X}}_z(f_{\psi}(z))=
\db{X}_z(e_{\bar{z}},\psi)= (e_{\bar{z}},\mb{X}\psi),
\end{equation}
where
\begin{equation}\label{sss}
\db{X}_{z}(f_{\psi}(z))=\left(P_{X}(z)+
\sum Q^i_X(z)\frac{\pa}{\pa  z_i}\right)f_{\psi}(z) .
\end{equation}
In \cite{last} we have advanced the hypothesis that for {\it CS-groups
the holomorphic functions $P$ and $Q$ in  {\em(\ref{sss})} are
polynomials},
 i.e. $\db{A} \subset \got{A}_1\subset \D_1$.

In this paper we present explicit formulas
for  (\ref{sss}) in case of the simplest example of a mixed group which
is a CS-group, the Jacobi group $G^J_1$. We start with the Jacobi algebra.

\section{The Jacobi algebra}\label{jac1}
The Heisenberg-Weyl   group  is the group with the
3-dimensional real  Lie algebra isomorphic to the Heisenberg  algebra
 \begin{equation}\label{nr0}\got{h}_1\equiv\got{g}_{HW}=
<is 1+xa^+-\bar{x}a>_{s\in\R ,x\in\C} ,\end{equation}
 where ${ a}^+$ (${ a}$) are  the boson creation
(respectively, annihilation)
operators which verify the CCR 
(\ref{baza1}).

Let us also consider the Lie algebra of the group $\text{SU}(1,1)$:
\begin{equation}\label{nr1}
\got{su}(1,1)=
<2i\theta K_0+yK_+-\bar{y}K_->_{\theta\in\R ,y\in\C} , \end{equation} 
where the generators $K_{0,+,-}$ verify the standard commutat
ion relations
(\ref{baza2}).
We  consider the matrix realization
\begin{equation}\label{nr2}
K_0=\frac{1}{2}\left(\begin{array}{cc}1 &0 \\0 &
-1\end{array}\right),~
K_+=i\left(\begin{array}{cc}0 &1 \\ 0 & 0 \end{array}\right)~,
K_-=i\left(\begin{array}{cc}0 &0 \\ 1 & 0 \end{array}\right) .
\end{equation}
 Now let us define the Jacobi algebra as the  the semi-direct sum
\begin{equation}\label{baza}
\got{g}^J_1:= \got{h}_1\rtimes \got{su}(1,1),
\end{equation}
where $\got{h}_1$ is an  ideal in $\got{g}^J_1$,
i.e. $[\got{h}_1,\got{g}^J_1]=\got{h}_1$, 
determined by the commutation relations:
\begin{subequations}\label{baza11}
\begin{eqnarray}
& & [a,a^+]=1\label{baza1}, \\
\label{baza2}
~& & \left[ K_0, K_{\pm}\right]=\pm K_{\pm}~,~ 
\left[ K_-,K_+ \right]=2K_0 , \\
\label{baza3}
& & \left[a,K_+\right]=a^+~,~\left[ K_-,a^+\right]=a, \\
& & \label{baza4}
\left[ K_+,a^+\right]=\left[ K_-,a \right]= 0 ,\\
\label{baza5}
& & \left[ K_0  ,~a^+\right]=\frac{1}{2}a^+, \left[ K_0,a\right]
=-\frac{1}{2}a .
\end{eqnarray}
\end{subequations}

\section{The differential action}\label{diff}
We shall suppose that we know the  derived representation $d\pi$ of 
the Lie algebra $\got{g}^J_1$ (\ref{baza})
of the Jacobi group $G^J_1$. We
associate to the generators $a, a^+$ of the HW group and to the
generators 
$K_{0,+,-}$ of the group  $\text{SU}(1,1)$ the operators $a, a^+$, respectively 
${\mb{K}}_{0,+,-}$, where 
$(a^+)^+=a,~ {\mb{K}}^+_0={\mb{K}}_0, {\mb{K}}^+_{\pm}={\mb{K}}_{\mp}$,
 and we impose to the cyclic vector $e_0$ to verify simultaneously
 the conditions
\begin{subequations}\label{cond}
\begin{eqnarray}
ae_0 & = & 0, \label{a0}\\
~ {\mb{K}}_-e_0 & = & 0, \label{kminus}\\
{\mb{K}}_0e_0 & = & k e_0;~ k>0, 2k=2,3,... .\label{k0}
\end{eqnarray}
\end{subequations}
We consider in  (\ref{k0}) the positive  discrete series
representations $D^+_k$ of $\text{SU}(1,1)$ (cf. \S 9 in \cite{bar47}).

 Perelomov's coherent state vectors   associated to the group $G^J_1$ with 
Lie algebra the Jacobi algebra (\ref{baza}), based on the manifold $M$:
\begin{subequations}\label{nmm}
\begin{eqnarray}\label{mm17}
M& := & H_1/\R\times \text{SU}(1,1)/\text{U}(1),\\
\label{mm}
M & = & \mc{D}^J_1:=\C\times\mc{D}_1, 
\end{eqnarray}
\end{subequations}
are defined as 
\begin{equation}\label{csu}
e_{z,w}:=e^{za^++w{\mb{K}}_+}e_0, ~z\in\C,~ |w|<1 .
\end{equation}
The general scheme (\ref{correspond}) associates to elements of the
Lie algebra $\got{g}$ differential operators:
 $X\in\got{g}\rightarrow\db{X}\in\D_1$. The space of
functions on which these operators act in the case of the Jacobi group
 will be made precise in \S \ref{sFs}.

The following lemma expresses the differential action of the generators
of the Jacobi  algebra as operators of the type $\got{A}_1$ in two
variables on $M$.
\begin{lemma}\label{mixt}The differential action of the generators
{\em (\ref{baza1})-(\ref{baza5})} of
the Jacobi algebra {\em (\ref{baza})} is given by the formulas:
\begin{subequations}\label{summa}
\begin{eqnarray}
& & \mb{a}=\frac{\pa}{\pa z};~\mb{a}^+=z+w\frac{\pa}{\pa z} ;\\
 & & \db{K}_-=\frac{\pa}{\pa w};~\db{K}_0=k+\frac{1}{2}z\frac{\pa}{\pa z}+
w\frac{\pa}{\pa w};\\
& & \db{K}_+=\frac{1}{2}z^2+2kw +zw\frac{\pa}{\pa z}+w^2\frac{\pa}{\pa
w} ,
\end{eqnarray}
\end{subequations}
where $z\in\C$, $|w|<1$.
\end{lemma}

{\it Proof}.
With the definition (\ref{csu}), we have the formal relations:
$$a^+e_{z,w}=\frac{\pa}
{\pa z}e_{z,w};~~{\mb{K}}_+e_{z,w}=\frac{\pa}{\pa w}e_{z,w}.$$ 
The proof is based on the general formula
\begin{equation}\label{for}
\Ad (\exp X) =\exp (\ad_X),
\end{equation}
valid  for Lie algebras $\got{g}$, which here  we write down explicitly as
\begin{equation}\label{bazap}
Ae^X=e^X(A-[X,A]+\frac{1}{2}[X,[X,A]]+\cdots ) ,
\end{equation}
and we take  $X=za^++w {\mb{K}}_+$ because of the definition (\ref{csu}).

1) Firstly we take in  (\ref{bazap})
$A=a$. Then $[X,A]=-z-wa^+$; $[X,[X,A]]=0 $, and
$$ae^X=e^X(a+z+wa^+) ;$$
$$ae^Xe_0=(z+w\frac{\pa}{\pa z})e^Xe_0 .$$

2) Now we take in (\ref{bazap})
 $A={\mb{K}}_0$. Then $[X,A]=-\frac{z}{2}a^+-w{\mb{K}}_+$;
$[X,[X,A]]=0 $, and
$${\mb{K}}_0e_{z,w}=(k+\frac{z}{2}\frac{\pa}{\pa z}+w\frac{\pa }{\pa
w})e_{z,w} .$$

3) Finally, we take in  (\ref{bazap})
$A={\mb{K}}_-$. We have  $[X,A]=-za-2w{\mb{K}}_0$, and
\begin{eqnarray*}
[X,[X,A]] & = & [za^++w {\mb{K}}_+,-za-2w{\mb{K}}_0]\\
 & = &
-z^2[a^+,a]-2zw[a^+,{\mb{K}}_0]-wz[{\mb{K}}_+,a]-2w^2[{\mb{K}}_+,{\mb{K}}_0]\\
 & = & z^2+2zwa^++2w^2{\mb{K}}_+ .
\end{eqnarray*}
Using  (\ref{bazap}), we have
\begin{eqnarray*}
Ae^Xe_0 & = &
e^X[{\mb{K}}_-+(za+2w{\mb{K}}_0)+\frac{1}{2}(z^2+2wza^++2w^2{\mb{K}}_+)]e_0\\
 & = & (2wk +\frac{z^2}{2}+wz\frac{\pa}{\pa z}+w^2
\frac{\pa}{\pa w})e^Xe_0  . 
\end{eqnarray*}

Now we do some general considerations. For
  $X\in\got{g}$, let $\mb{X}.e_{z}:=\mb{X}_ze_z$. Then
$\mb{X}.e_{\bar{z}}={\mb{X}}_{\bar{z}}e_{\bar{z}}$.

 But
$(e_{\bar{z}},
{\mb{X}}.e_{\bar{z}'})=({\mb{X}}^+.e_{\bar{z}},e_{\bar{z}'})$ and
finally, with equation (\ref{car1}), we have
\begin{equation}\label{bidiff}
\db{X}_{\bar{z}'}(e_{\bar{z}},e_{\bar{z}'})=\db{X}^+_{z}
(e_{\bar{z}},e_{\bar{z}'})
\end{equation}
With observation (\ref{bidiff}) and the previous calculation,
 Lemma \ref{mixt} is proved.   \hfill$\gata$

\begin{Comment}
We illustrate {\em   (\ref{bidiff})} for $\mb{X}=\mb{a}$. Then it can be
checked up that
$$(\bar{z}'+\bar{w}'\frac{\pa}{\pa \bar{z}'})
(e_{\bar{z},\bar{w}},e_{\bar{z}',\bar{w}'})=\frac{\pa}{\pa z}
(e_{\bar{z},\bar{w}},e_{\bar{z}',\bar{w}'})=
\frac{\bar{z}'+z\bar{w}'}{1-w\bar{w}'} 
  (e_{\bar{z},\bar{w}},e_{\bar{z}',\bar{w}'}) , $$
where the kernel has the expression {\em (\ref{KHK})} calculated  below.
\end{Comment}

\section{The reproducing kernel}\label{repK}

Now we calculate the reproducing kernel $K$ on the base manifold
$ M=\mc{D}^J_1$ as the scalar product of two Perelomov's
CS-vectors  (\ref{csu}), taking into account the conditions
(\ref{cond}) and the orthonormality of the basis of the Hilbert spaces
associated with the factors of the Jacobi group.
\begin{lemma}\label{lema5}
Let $K=K(\bar{z},\bar{w},z,w)$, where $z\in\C, ~w\in\C, ~|w|<1$,
\begin{equation}\label{KHK1}
K:=(e_0,e^{\bar{z}a+\bar{w}{\mb{K}}_-}e^{za^++w{\mb{K}}_+}e_0) .
\end{equation}
Then the reproducing kernel is
\begin{equation}\label{hot}
K=
(1-w\bar{w})^{-2k}\exp{\frac{2z\bar{z}+z^2\bar{w}+\bar{z}^2w}{2(1-w\bar{w})}}. 
\end{equation}
More generally, the kernel $K:\mc{D}^J_1\times \bar{\mc{D}}^J_1\rightarrow \C$ is: 
\begin{equation}\label{KHK}
K(z,w;\bar{z}',\bar{w}'):=(e_{\bar{z},\bar{w}},e_{\bar{z}',\bar{w}'})
=(1-{w}\bar{w}')^{-2k}\exp{\frac{2\bar{z}'{z}+z^2\bar{w}'
+\bar{z}'^2w}{2(1-{w}\bar{w}')}} .
\end{equation}
\end{lemma}

{\it Proof}.
We introduce the auxiliary  operators:
\begin{subequations}
\begin{eqnarray}
\label{kkk1}{\mb{K}}_+ & = &\frac{1}{2}(a^+)^2+{\mb{K}}'_+ ,\\
\label{kkk2}{\mb{K}}_- & = &\frac{1}{2}a^2+{\mb{K}}'_- ,\\
{\mb{K}}_0 & = &\frac{1}{2}(a^+a+\frac{1}{2})+{\mb{K}}'_0 ,
\end{eqnarray}
\end{subequations}
which have the properties
\begin{subequations}
\begin{eqnarray}
{\mb{K}}'_-e_0 & = & 0 ,\\
{\mb{K}}'_0e_0& = & k'e_0; ~k=k'+\frac{1}{4};
\end{eqnarray}
\end{subequations}
\begin{subequations}
\begin{eqnarray}
\label{kkk3} & & \left[{\mb{K}}'_{\sigma},a\right]=\left[
{\mb{K}}'_{\sigma},a^+\right]=0,~\sigma =\pm ,0 ,\\ 
\label{kkk4} & & \left[ {\mb{K}}'_0,{\mb{K}}'_{\pm}\right]=\pm {\mb{K}}'_{\pm};~
 \left[ {\mb{K}}'_-,{\mb{K}}'_+\right]=2{\mb{K}}'_0 .
\end{eqnarray}
\end{subequations}

 Using the fact that $e_{k,k+m}$ {\it is an orthonormal
system} (see also \S \ref{su1} and the Appendix), where
\begin{equation}\label{ort1}
e_{k,k+m}:=a_{km}({\mb{K}}_+)^me_{k,k};~~a^2_{km}=\frac{\Gamma
(2k)}{m!\Gamma (m+2k)},
\end{equation}
 the relation (see e.g. equation  1.110  in \cite{grad})
\begin{equation}\label{rajik}
(1-x)^{-q}=\sum_{m=0}^{\infty}\frac{x^m}{m!}\frac{\Gamma
(q+m)}{\Gamma (q)},
\end{equation}
and  the orthonormality of the $n$-particle states (see also \S
\ref{heil} and the Appendix):
\begin{equation}\label{ort2}
|n>=(n!)^{-\frac{1}{2}}(a^+)^{n}|0>;~~<n',n>=\delta_{nn'} ,
\end{equation} 
it is proved the relation
\begin{equation}
\label{una1}
(e_0,e^{\bar{w}{\mb{K}}'_-}e^{w'{\mb{K}}'_+}e_0)=(1-w'\bar{w})^{-2k'} .
\end{equation}

We introduce the notation
\begin{equation}\label{ee}
E=E(z,w):=  e^{za^++\frac{w}{2}(a^+)^2}=
\sum_{p,q\ge 0}\frac{z^p}{p!}\frac{{(\frac{w}{2}})^q}{q!}(a^+)^{p+2q}.
\end{equation}
With the change of variable: 
$n:=p+2q$, i.e. $ p=n-2q$, equation  (\ref{ee}) becomes
$$E=\sum_{n\ge 0}\sum_{q=0}^{\left[ \frac{n}{2}\right]}
\frac{z^{n-2q}}{(n-2q)!q!}{\left(\frac{w}{2}\right)}^q(a^+)^n .$$
Recalling that  the Hermite polynomials can be represented as
(cf. equation 10.13.9 in \cite{bate})
\begin{equation}\label{hermite} 
H_n(x)=n!\sum_{m=0}^{\left[\frac{n}{2}\right]}
\frac{(-1)^m(2x)^{n-2m}}{m!(n-2m)!} ,
\end{equation}
the expression (\ref{ee}) becomes 
\begin{equation}\label{cucu1}
E(z,w)=
\sum_{n\ge 0}\frac{i^{-n}}{n!}\left(\frac{w}{2}\right)^{\frac{n}{2}}
H_n\left(\frac{iz}{\sqrt{2w}}\right)(a^+)^{n} .
\end{equation}
Then
$$K:=
K(\bar{z},\bar{w};z',w')=(e_{z,w},e_{z',w'})=
(e_0,e^{\bar{z}a+\bar{w}{\mb{K}}_-}
e^{z'a^++w'{\mb{K}}_+}e_0) .$$
But due to  equations (\ref{kkk1}), (\ref{kkk2}),
 $K$ can be written down as
$$K=(e_0,e^{\bar{z}a+\frac{\bar{w}}{2}a^2}e^{\bar{w}{\mb{K}}'_-}
e^{w'{\mb{K}}'_+}e^{z'a^++\frac{w'}{2}(a^+)^2}e_0) .$$
Let the notation
$$F:=
F(\bar{z}'\bar{w}';z,w)= (e_0,E^+(\bar{z}',\bar{w}')E(z,w)e_0) .$$
Because of the 
 orthonormality  relation (\ref{ort2}),
$(e_0,a^{n'}(a^+)^ne_0)=n!\delta_{nn'}$, we get:   
$$ F
=\sum \frac{1}{n!}
\left(\frac{\bar{w}'w}{4}\right)^{\frac{n}{2}}
H_n(-i\frac{\bar{z}'}{\sqrt{2\bar{w}'}})H_n(i\frac{z}{\sqrt{2w}}).$$
We use the summation relation of the Hermite polynomials
 (Mehler formula, cf. equation 10.13.22 in  \cite{bate})
\begin{equation}\label{bat}
\sum_{n=o}^{\infty}\frac{(\frac{s}{2})^n}{n!}H_n(x)H_n(y)=
\frac{1}{\sqrt{1-s^2}}\exp{\frac{2xys-(x^2+y^2)s^2}{1-s^2}},~ |s|<1 , 
\end{equation}
where
$$x=-i\frac{\bar{z}'}{\sqrt{2\bar{w}'}};~
y=i\frac{z}{\sqrt{2w}};~ s=(\bar{w}'w)^{1/2},$$
and we get
$$F
=\frac{1}{{(1-\bar{w}'w)}^{1/2}}
\exp\frac{2\bar{z}'z+\bar{z}'^2w+z^2\bar{w}'}{2(1-\bar{w}'w)} .$$
Recalling  (\ref{una1}), we have
 $$K=(1-\bar{w}'w)^{-2k'}F, $$
and finally:
$$(e_{\bar{z},\bar{w}},e_{\bar{z}',\bar{w}'})=
(1-w\bar{w}')^{-2k}\exp \frac{2z\bar{z}'+z^2\bar{w}'
+\bar{z}'^2w}{2(1-w\bar{w}')}.  $$
{\mbox{~~~}} \hfill $\gata$

\section{The group action on the base manifold}\label{gra}

We start this section  recalling in \S \ref{unul}
  some useful relations for
 representations of the
groups $H_1$ and $\text{SU}(1,1)$. Then  we obtain  formulas (\ref{schimb0}),
(\ref{schimb}) for the change of
order of the action of these groups.

\subsection{Formulas for the Heisenberg-Weyl group $H_1$
 and $\text{SU}(1,1)$}\label{unul}
Let us recall some relations for the displacement operator:
\begin{equation}\label{deplasare}
D(\alpha ):=\exp (\alpha a^+-\bar{\alpha}a)=\exp(-\frac{1}{2}|\alpha
|^2)  \exp (\alpha a^+)\exp(-\bar{\alpha}a),
\end{equation}
\begin{equation}\label{thetah}
D(\alpha_2)D(\alpha_1)=e^{i\theta_h(\alpha_2,\alpha_1)}
D(\alpha_2+\alpha_1) , 
~\theta_h(\alpha_2,\alpha_1):=\Im (\alpha_2\bar{\alpha_1}) .
\end{equation}

Let us denote by $S$ the $D^k_+$ representation of the group $\text{SU}(1,1)$
and  let us introduce the notation $\underline{S}(z)=S(w)$, where
$w$ and $z$, $w\in\C,~
|w|<1$,  $z\in\C$, 
are related by  (\ref{u3}), (\ref{u4}). We have the relations:
\begin{subequations}
\begin{eqnarray}
\underline{S}(z) & := &\exp (z{\mb{K}}_+-\bar{z}{\mb{K}}_-), ~z\in\C
;\label{u1} \\
S(w) & = &  \exp (w{\mb{K}}_+)\exp (\eta
{\mb{K}}_0)\exp(-\bar{w}{\mb{K}}_-);\label{u2}\\
w & = & w(z)  = \frac{z}{|z|}\tanh \,(|z|)~, w\in\C, |w|<1;\label{u3}\\
z & = & z(w)  =  \frac{w}{|w|}
{\text{arctanh}} \,(|w|) = \frac{w}{2|w|}\log\frac{1+|w|}{1-|w|};\label{u4}\\
\eta & =& \log (1-w\bar{w})= -2\log (\cosh \,(|z|)) . \label{u5}
\end{eqnarray}
\end{subequations}
Let us consider an element $g\in \text{SU}(1,1)$, 
\begin{equation}\label{dg}
g= \left( \begin{array}{cc}a & b\\ \bar{b} &
\bar{a}\end{array}\right),~\text{where} ~|a|^2-|b|^2=1.
\end{equation}
\begin{Remark}The following relations hold:
\begin{equation}\label{r1}
\underline{S}(z)e_0=(1-|w|^2)^k e_{0,w} ,
\end{equation}
\begin{equation}\label{r2}
e_g:= S (g) e_0
= \bar{a}^{-2k}e_{0, w=-i\frac{b}{\bar{a}}}=
\left(\frac{a}{\bar{a}}\right)^{k}\underline{S}(z)e_0 ,
\end{equation}
\begin{equation}\label{r3}
S(g)e_{0,w}=(\bar{a}+\bar{b}w)^{-2k}e_{0,g\cdot w}, 
\end{equation}
where $w\in\C,~ |w|<1$ and  $z\in\C$
 in {\em  (\ref{r2})} are related by equations  {\em   (\ref{u3})},
{\em (\ref{u4})},  and
the linear-fractional action of the group $\rm{SU}(1,1)$ on
the unit disk $\mc{D}_1$ in {\em   (\ref{r3})} is
\begin{equation}\label{r4}
g\cdot w 
=\frac{a \, w+ b}{\bar{b}\, w +\bar{a}}.
\end{equation}
\end{Remark}

We recall also the following property, which is a particular case of a
more general result proved   in \cite{cluj}:
\begin{Remark}If $\underline{S}(z)$ is defined by {\em{(\ref{u1})}}, then: 
\begin{subequations}
\begin{eqnarray}
\underline{S}(z_2)\underline{S}(z_1) & = & \underline{S}(z_3)e^{i\theta_s
{\mb{K}}_0}\label{unuu1};\\
w_3 & = & \frac{w_1+w_2}{1+\bar{w}_2w_1};\label{unuu2}\\
e^{i\theta_s} & = & \frac{1+w_2\bar{w}_1}{1+w_1\bar{w}_2},
\label{thetas}
\end{eqnarray}
where  $w_i$ and $z_i$, $i=1,2,3$, in equation {\em  (\ref{unuu2})}
 are related by the relations {\em (\ref{u3}), (\ref{u4})}.
\end{subequations}
\end{Remark}
\begin{Comment}
Note that when $z_1,z_2\in\R$, then {\em  (\ref{unuu1})} expresses just 
the additivity of the ``rapidities'', 
$$\underline{S}(z_2)\underline{S}(z_1)  =  \underline{S}(z_2+z_1),$$
while {\em   (\ref{unuu2})} becomes just the Lorentz composition of
velocities in special relativity: $$w_3= \frac{w_1+w_2}{1+w_2w_1} .$$
\end{Comment}

\subsection{Holstein-Primakoff-Bogoliubov-type equations}\label{hpb}

We recall the {\it Holstein-Pri\-makoff-Bogoliubov  equations}
\cite{hol},\cite{bog} (see also \cite{nieto}), a consequence
of the equation (\ref{for}) and of the fact that the Heisenberg algebra
is an ideal in the Jacobi algebra (\ref{baza}), as expressed in
 (\ref{baza3})-(\ref{baza5}): 
\begin{subequations}\label{hol}
\begin{eqnarray}
\underline{S}^{-1}(z)\, a \, \underline{S}(z) & = &
 \cosh(|z|)\, a + \frac{z}{|z|}\sinh (|z|)\,a^+,\\
\underline{S}^{-1}(z)\, a^+ \, \underline{S}(z) & = &
\cosh (|z|)\, a^+ + \frac{\bar{z}}{|z|}\sinh (|z|)\, a  ,
\end{eqnarray}
\end{subequations}
and the CCR are still fulfilled in the new creation and annihilation operators.

Let us introduce the notation: 
\begin{equation}\label{notatie}
\tilde{A}=\left(\begin{array}{c}
A\\ \bar{A}\end{array}\right);~~
\mc{D}={\mc{ D}}(z) = 
\left( \begin{array}{cc} M & N \\ P& Q\end{array}  \right),
\end{equation}
where
\begin{equation}\label{notatie1} M  =  \cosh \,(|z|) ;~
N  =  \frac{z}{|z|}\sinh\,(|z|) ;~
P  =  \bar{N} ;~
Q  =   M .
\end{equation}
Note that 
\begin{equation}
{\mc{ D}}(z) = e^{X},~\text{where} ~
 X:=\left(\begin{array}{cc} 0 & z\\ \bar{z} &
0\end{array}\right)  .
\end{equation}
\begin{Remark}
With the notation {\em (\ref{notatie}), (\ref{notatie1})}, equations
    {\em(\ref{hol})} become:
$$\underline{S}^{-1}(z)\tilde{a}\underline{S}(z)={\mc{D}}(z)\tilde{a}.$$
\end{Remark}

Using formula (\ref{for}), we obtain, as a consequence that the  HW
group 
is a normal subgroup of the Jacobi group, the relations (\ref{schimb0}),
(\ref{schimb}) (or (\ref{schimb2})), which allow to
interchange the order of the representations of the groups $\text{SU}(1,1)$
and HW:
\begin{Remark}\label{R44} If $D$ and $\underline{S}(z)$ are defined by
{\em{(\ref{deplasare})}}, respectively {\em{(\ref{u1})}}, then 
\begin{equation}\label{schimb0}
D(\alpha )\underline{S}(z)= \underline{S}(z) D(\beta ),
\end{equation}
where
\begin{subequations}\label{schimb}
\begin{eqnarray}
\label{schimb33} \beta =\alpha \cosh\,
(|z|)-\bar{\alpha}\frac{z}{|z|}\sinh \,(|z|) & ; & \alpha = \beta \cosh \,(|z|)
+ \bar{\beta}\frac{z}{|z|}\sinh \,(|z|);\\
\label{schimb34}\beta = 
\frac{\alpha - \bar{\alpha} w }{ (1- |w|^2)^{1/2} } & ; &
\alpha = \frac{\beta +\bar{\beta} w}{(1- |w|^2)^{1/2} }.
\end{eqnarray}
\end{subequations}
With the convention {\em (\ref{notatie})}, equation  {\em (\ref{schimb33})}
 can be written down as:
\begin{equation}\label{schimb2}
\tilde{\beta}={\mc{D}}(-z)\tilde{\alpha};~
 \tilde{\alpha}={\mc{D}}(z)\tilde{\beta}.
\end{equation}
\end{Remark}

Let us introduce  the notation
\begin{equation}\label{k00}
\underline{S}(z,\theta ):= \exp (2i\theta
{\mb{K}}_0+z{\mb{K}}_+-\bar{z}{\mb{K}}_-). 
\end{equation}

Using  (\ref{for}),    more general formulas than 
Holstein-Primakoff-Bogoliubov  equations (\ref{hol}) can be proved, namely:
\begin{subequations}\label{hol1}
\begin{eqnarray}
\underline{S}(z,\theta )^{-1}(z)\,a\,\underline{S}(z,\theta ) & = &
(\text{cs}\,(x)+ i \theta \,\frac{\text{si}\,(x)}{x})\, a  + 
z\frac{\text{si} \,(x)}{x} \,  a^+,\\
\underline{S}(z,\theta )^{-1}(z)\,a^+ \, \underline{S}(z, \theta ) & = &
(\text{cs} \,(x)-i\theta\,\frac{\text{si} \,(x)}{x}) \, a^+
+\bar{z}\frac{\text{si} \,(x)}{x} \, a ,
\end{eqnarray}
\end{subequations}
where 
\begin{equation}
{\text{cs}}\, (x) := \left\{
\begin{array}{ll} \cosh \,(x), & {\text{if~}}
\lambda =x^2>0 , \\
\cos \,(x), &  {\text{if~}}
\lambda = -x^2<0, 
 \end{array}
\right. ; \lambda := |z|^2-\theta^2, 
\end{equation}
and similarly for $\text{si}\,(x)$.

 Let us
consider $X\in\got{su}(1,1)$,
\begin{equation}
X=\left( \begin{array}{cc} i\theta & z \\ \bar{z} & -i\theta 
\end{array} \right)~,~ \theta\in\R ,~ z\in\C . 
\end{equation}
Then  $g=e^X\in \text{SU}(1,1)$ is an element of the form (\ref{dg}), where
\begin{equation}\label{zzz}
  a = \text{cs}\, (x) + i\theta\, \frac{\text{si}\, (x)}{x}, ~
b = z \frac{\text{si} \, (x) }{x} .
\end{equation}
If $g= \left(\begin{array}{ll}\alpha & \beta \\ \bar{\beta} &\bar{\alpha}
\end{array}\right)\in \text{SU}(1,1)$,
  then equations (\ref{hol1}) can be written down as 
\begin{subequations}\label{hol2}
\begin{eqnarray}
S^{-1}(g)\,a\,S(g) & = & \alpha\, a +\beta \, a^+, \\
S^{-1}(g)\,a^+\,S(g) & = & \bar{\beta}\, a +\bar{\alpha} \, a^+ ,
\end{eqnarray}
\end{subequations}
and we have  the following
(generalized Holstein-Primakoff-Bogoliubov) equations:
\begin{Remark} If $S$ denotes the representation of $\rm{SU}(1,1)$, with
the convention {\em{(\ref{notatie})}},  we have 
\begin{equation}\label{hol3}
S^{-1}(g)\,\tilde{a}\, S (g)= g\cdot \tilde{a}.
\end{equation}
\end{Remark}

Applying again  formula (\ref{for}), we obtain a
 more general formula than  (\ref{schimb0}), namely:
\begin{equation}\label{maimult}
\underline{S}(z,\theta)D (\alpha )\underline{S}(z,\theta )^{-1}= D(\alpha_1),
\end{equation}
where 
\begin{equation}\label{alpha1}
\alpha_1= \alpha_1(z,\alpha , \theta )
=  \alpha\, \text{cs} \,(x)+
(i\theta \alpha +z\bar{\alpha})\,\frac{\text{si}\, (x)}{x} .
\end{equation} 

Written down in the form similar to  (\ref{schimb0}), 
equation (\ref{maimult}) reads
\begin{equation}\label{maimult1}
D(\alpha ) \underline{S}(z,\theta )=\underline{S}(z,\theta )
D(\beta_1) , 
\end{equation}
where $\beta_1 = \beta_1(z,\alpha , \theta )=\alpha_1(z,-\alpha ,
-\theta )$, i.e.
\begin{equation}\label{schimb3}
\beta_1=\alpha\,{\text{cs}}\,(x)-
(i\theta\alpha  +z\bar{\alpha})\,\frac{\text{si}\,(x)}{x}~, \alpha =
\beta_1\,\text{cs}\,(x)+(i\theta\beta_1 
+z\bar{\beta}_1)\,\frac{\text{si}\,(x)}{x}.
\end{equation}
Note that if $\theta =0$, then $\underline{S}(z,\theta )=
\underline{S}(z)$ and $\beta_1$ in  (\ref{schimb3}) becomes
$\beta_1=\beta$ with $\beta$ given by (\ref{schimb}).

We also underline that if $z=0$ in  (\ref{maimult}), then
 (\ref{alpha1}) becomes just
$$\alpha_1= \alpha (\cos \,(|\theta|)+i\theta 
\frac{\sin \,(|\theta |)}{ |\theta |}). $$

Summarizing, we rewrite now  equation
  (\ref{maimult}) in the following  useful form:
\begin{Remark}\label{remark6}
In the matrix realization {\em (\ref{nr2})}, equation  
{\em(\ref{maimult})} can be written down as
\begin{equation}\label{maimult3}
S (g) D(\alpha ) S^{-1}(g) = D (\alpha_g),
\end{equation}
where {\em (\ref{alpha1})} has the expression of the natural action of
$\rm{SU}(1,1)\times \C\rightarrow \C$: $g\cdot\tilde{\alpha}:=\alpha_g$, 
\begin{equation}\label{alpha2}
\alpha_g = a\,\alpha + b\,\bar{\alpha}, 
\end{equation}
and $a$, $b$ have the expression {\em (\ref{zzz})}.
\end{Remark}

Let us remark that 
the commutation relations
(\ref{baza3})-(\ref{baza5}) between the generators of the groups
$\text{SU}(1,1)$ and  HW   were chosen in such a way that the
action of the group $\text{SU}(1,1)$ on the complex plane $M\approx
\C=H_1/\R$ be the natural one,
cf. Remark \ref{remark6}. Such a choice of the action of the
group $\text{SU}(1,1)$ on the group $H_1$, a normal subgroup of the Jacobi
group $G^J_1$, was inspired from the squeezed states
  in Quantum
Optics (cf. e.g. \cite{nieto}). If we had started from the natural action of
 $\text{SU}(1,1)$ on $\C$
given in Remark \ref{remark6}, then the commutation relations
(\ref{baza3})-(\ref{baza5}) would had  followed taking the derivatives in
 (\ref{hol3}) realized as  (\ref{hol1}) using the development (\ref{zzz}).

Now we consider the product of two representations $D$ and $S$ and
apply  Remark   \ref{R44}: 
\begin{equation}\label{alp1}
D(\alpha_2)\underline{S}(z_2)D(\alpha_1)\underline{S}(z_1)=
D(\alpha_2)D(\alpha )\underline{S}(z_2)\underline{S}(z_1),\end{equation}
where \begin{equation}\label{alp}
\alpha =\alpha_1\cosh \,(|z_2|)+\bar{\alpha}_1\frac{z_2}{|z_2|}\sinh\,
(|z_2|),
\end{equation}
or $$ \tilde{\alpha} = {\mc{D}}(z)\tilde{\alpha}_1.$$

Equations (\ref{alp1}) and (\ref{alp}) allow to determine
\begin{Remark}\label{actiunea11} The action:
$(\alpha_2,z_2)\times 
(\alpha_1, w_1)=(A,w)$, where $z_2,\alpha_{1,2}, A\in\C$, $ w,
w_1\in\mc{D}_1$
and the variables of type $w$ and $z$ are related by equations
{\em  (\ref{u3}), (\ref{u4})},  can be expressed as: 
\begin{subequations}\label{actiunea}
\begin{eqnarray}
A & = & \alpha_2 +\alpha_1 \cosh |z_2|
+\bar{\alpha}_1\frac{z_2}{|z_2|}\sinh |z_2| = \alpha_2 +
\frac{\alpha_1+\bar{\alpha}_1 w_2}{ (1-|w_2|^2)^{1/2}},\\
w & = & \frac{\cosh|z_2|w_1 +\frac{z_2}{|z_2|}\sinh |z_2|}
{\frac{\bar{z}_2}{|z_2|}\sinh |z_2|w_1 +\cosh |z_2|} =
\frac{w_1+w_2}{1+ w_1\bar{w}_2}.
\end{eqnarray}
\end{subequations}
 Equations {\em (\ref{actiunea})} express the action $(\alpha_2,w_2)\times
(\alpha_1, w_1)= (\alpha_2 + w_2\circ\alpha_1, w_2\circ w_1)$,
$\alpha_{1,2}\in\C,
w_{1,2}\in \mc{D}_1$.  
{\em  (\ref{actiunea})} can be written down as:
\begin{subequations}
\begin{eqnarray}
\tilde{A} & = & \tilde{\alpha}_2+ {\mc{D}}\tilde{\alpha}_1,\\
w & = & \frac{M w_1+ N}{P w_1 + Q} .
\end{eqnarray}
\end{subequations}
\end{Remark}
Let us introduce the normalized vectors:
\begin{equation}\label{psi2}
\Psi_{\alpha, w}:=D(\alpha ) S(w) e_0;~\alpha\in\C, ~w\in\C, |w|<1. 
\end{equation}
As a consequence of  (\ref{alp1}), we have:
\begin{Remark} The product of the representations $D$
and $\underline{S}$ acts on the CS-vector {\em{(\ref{psi2})}} with the effect:
\begin{equation}
D(\alpha_2)\underline{S}(z_2)\Psi_{\alpha_1,w_1}=J\Psi_{A,w},
~\text{where}~J= e^{i(\theta_h(\alpha_2,\alpha ) + k\theta_s )}.
\end{equation}
Above $(A,w)$ are given by {\em Remark \ref{actiunea11}},
 $\theta_h(\alpha_2,\alpha)$ is given by {\em  (\ref{thetah})} with
$\alpha $ given by {\em  (\ref{alp})},
while $\theta_s$ is given by {\em  (\ref{thetas})} and  the dependence
$w_2=w_2(z_2)$ is
given by equation {\em  (\ref{u3})}.
\end{Remark}
Note also  the following important  property   (\ref{stolo}),
 well known in the Quantum Optics of
squeezed states  (see e.g.  equation (20) p. 3219 in \cite{stol}):

\begin{Comment}\label{STO} {\it The action of the HW
group on the (``squeezed'') state vector}
$$\underline{\Psi}_{z,\alpha}=\underline{S}(z) D({\alpha})e_0$$ 
 {\it modifies only  the part of the
HW group. More precisely, we have}
\begin{equation}\label{stolo}
D(\beta )\underline{\Psi}{z,\alpha }= e^{i\eta}
\underline{\Psi}{z,\alpha +\gamma},
 ~\text{where}~\eta =\Im (\gamma\bar{\alpha}), 
\end{equation}
and 
\begin{equation}\label{stolo1}
\gamma =\beta \cosh (|z|)-\bar{\beta}\frac{z}{|z|}\sinh (|z|), ~
\text{or}~
\tilde{\gamma}={\mc{D}}(-z)\tilde{\beta}.
\end{equation}
\end{Comment}

 Indeed, we apply formula (\ref{schimb0}):
$$D(\beta )\underline{S}(z)D(\alpha ) = \underline{S}(z)D(\gamma)
D(\alpha ) ,$$
where $\gamma$ has the expression (\ref{stolo1}). 

Then  (\ref{stolo}) follows. \hfill $\gata$

\subsection{The  action of the Jacobi group}\label{actg}

Now we find  a relation  between the (normalized) vector
(\ref{psi2}) and the (unnormalized) Perelomov's CS-vector (\ref{csu}),  
 which will be important in the proof of Proposition
\ref{mm1}, our   main result of this section.
\begin{lemma}\label{lema6} The vectors {\em (\ref{psi2}), (\ref{csu})},
i.e.
$$\Psi_{\alpha, w}:=D(\alpha ) S(w) e_0;~ e_{z,w'}:=\exp
(za^++w'{\mb{K}}_+)e_0 . $$
are related by the relation
\begin{equation}\label{csv}
\Psi_{\alpha, w} = (1-w\bar{w})^k
\exp (-\frac{\bar{\alpha}}{2}z)e_{z,w}  ,
\end{equation}
where $z=\alpha-w\bar{\alpha}$.
\end{lemma}

{\it Proof}.
Due to  (\ref{u1}), (\ref{u2}), (\ref{kminus}) and (\ref{k0}),
we have the relations
\begin{eqnarray*}
S(w)e_0 & = &\exp (w{\mb{K}}_+)
\exp (\eta {\mb{K}}_0)\exp (-\bar{w}{\mb{K}}_-)e_0\\
  & = & \exp (w{\mb{K}}_+)\exp(k\ln (1-w\bar{w}))\e_0\\
& = & (1-w\bar{w})^k\exp (w{\mb{K}}_+)e_0 ,
\end{eqnarray*}
which is also a proof of  (\ref{r1}).

We obtain  successively
\begin{eqnarray*}
\Psi_{\alpha, w}& = &(1-w\bar{w})^k D(\alpha ) \exp (w{\mb{K}}_+)e_0\\
 &= &(1-w\bar{w})^k\exp(-\frac{1}{2}|\alpha
|^2)
\exp (\alpha a^+)\exp(-\bar{\alpha}a)\exp (w{\mb{K}}_+)e_0\\
 &  = &(1-w\bar{w})^k\exp(-\frac{1}{2}|\alpha
|^2)\exp (\alpha a^+)\exp(-\bar{\alpha}a)\exp
(w{\mb{K}}_+)\exp(\bar{\alpha }a)\exp(-\bar{\alpha }a) e_0\\
 & = &(1-w\bar{w})^k \exp(-\frac{1}{2}|\alpha
|^2)\exp (\alpha a^+) E e_0 ,
\end{eqnarray*}
where here
\begin{equation}\label{eeee}
E:=\exp(-\bar{\alpha}a)\exp
(w{\mb{K}}_+)\exp(\bar{\alpha }a) .
\end{equation}
As a consequence of  (\ref{for}) 
$$\exp (Z)\exp (X) \exp (-Z)=\exp (X+[Z,X]+\frac{1}{2}[Z,[Z,X]]+\cdots ),$$
where, if we take
$Z=-\bar{\alpha}a;~ X =w{\mb{K}}_+$, then 
$$ [Z,X]= -\bar{\alpha}wa^+; [Z,[Z,X]]=\bar{\alpha}^2w .$$
We find for $E$ defined by (\ref{eeee}) the value
$$E= \exp w({\mb{K}}_+-\bar{\alpha}a^++\frac{\bar{\alpha}^2}{2}) ,$$
and finally
$$\Psi_{\alpha, w} = 
\exp (-\frac{1}{2}|\alpha
|^2)\exp(w\frac{\bar{\alpha}^2}{2})(1-w\bar{w})^k
e_{\alpha -w\bar{\alpha},w}, $$
i.e.  (\ref{csv}).\hfill $\gata$
\begin{Comment}Starting from {\em  (\ref{csv})}, 
we reobtain the expression   {\em (\ref{hot})} of the reproducing kernel
$K$.
\end{Comment}
Indeed,   the normalization
$(\Psi_{\alpha ,w},\Psi_{\alpha , w})=1$ implies  that
\begin{equation}\label{are}
(e_{\alpha- w\bar{\alpha}, w},e_{\alpha - w\bar{\alpha}, w})=\exp
(|\alpha|^2-\frac{w}{2}\bar{\alpha}^2-c.c.)(1-w\bar{w})^{-2k}.
\end{equation}
With  the notation:
$\alpha -w \bar{\alpha} = z,$ we have 
$$ \alpha = 
\frac{z+\bar{z}w}{1-\bar{w}w}, $$
and then (\ref{are}) can be rewritten as
$$(e_{z,w},e_{z,w})=(1-\bar{w}w)^{-2k}
\exp\frac{2z\bar{z}+w\bar{z}^2+\bar{w}z^2}{2(1-w\bar{w})},$$ 
i.e. we get another proof of  (\ref{hot}). \hfill $\gata$

From the following proposition we can  see the holomorphic 
action of the Jacobi  group
\begin{equation}\label{jac}
 G^J_1:=H_1\rtimes \text{SU}(1,1), 
\end{equation}
on the manifold $\mc{D}^J_1$  (\ref{mm}):

\begin{Proposition}\label{mm1}
Let us consider the action $S(g)D(\alpha )e_{z,w}$, where $g\in
\rm{SU}(1,1)$ has the form {\em (\ref{dg})},  $D(\alpha )$ is given by {\em
 (\ref{deplasare})}, and the
coherent state vector is defined in {\em  (\ref{csu})}. Then we have the
formula {\em (\ref{xx})} and the relations {\em (\ref{xxx}), 
(\ref{x1})-(\ref{x3})}:
\begin{equation}\label{xx}
S(g)D(\alpha )e_{z,w}=\lambda e_{z_1,w_1}, ~ \lambda = \lambda
(g,\alpha; z,w) , \end{equation}
\begin{equation}\label{xxx}
z_1=\frac{\alpha-\bar{\alpha}w+z}{\bar{b}w+\bar{a}}; ~ w_1=g\cdot
w=\frac{a w+ b}{\bar{b}w+\bar{a}},\end{equation}
\begin{equation}\label{x1}
\lambda= (\bar{a}+\bar{b}w)^{-2k}
\exp (\frac{z}{2}\bar{\alpha}_0
-\frac{z_1}{2}\bar{\alpha}_2)\exp i\theta_h(\alpha_1 ,\alpha_0),\end{equation}
\begin{equation}\label{x2}\alpha_0 =
\frac{z+\bar{z}w}{1-w\bar{w}},\end{equation} 
\begin{equation}\label{x3}
\alpha_2= (\alpha +\alpha_0)a +(\bar{\alpha}+\bar{\alpha}_0)b. \end{equation}
\end{Proposition}

\begin{corollary}The action of the 6-dimensional Jacobi group {\em
(\ref{jac})} on the 
4-dimensional manifold {\em (\ref{mm})}, where
$\mc{D}_1=\rm{SU}(1,1)/\rm{U}(1)$, is given by 
equations {\em (\ref{xx}), (\ref{xxx})}. The composition law in $G^J_1$ is
\begin{equation}\label{compositie}
(g_1,\alpha_1,t_1)\circ (g_2,\alpha_2, t_2)= (g_1\circ g_2,
g_2^{-1}\cdot\tilde{\alpha}_1+\alpha_2, t_1+ t_2 +\Im
(g^{-1}_2\cdot\alpha_1\bar{\alpha}_2)),
\end{equation}
where $g\cdot\tilde{\alpha} :=\alpha_g$
 is given by {\em  (\ref{alpha2})}, and if
$g$ has the form given by {\em  (\ref{dg})},  then
$g^{-1}\cdot\tilde{\alpha} =\alpha_{g^{-1}} = \bar{a}\alpha -b\bar{\alpha}$. 
\end{corollary}

{\it Proof of {\em  Proposition \ref{mm1}}}.
With  Lemma \ref{lema6}, we have
$e_{z,w}=\lambda_1\Psi_{\alpha_0, w}$, where $\alpha_0$ is given by
 (\ref{x2}) and $\lambda_1=\exp
(\frac{z}{2}\bar{\alpha}_0)(1-|w|^2)^{-k}$.
 Then $I:= S(g)D(\alpha)e_{z,w}$ becomes successively
\begin{eqnarray*}
 I & = & \lambda_1 S(g) D(\alpha )\Psi_{\alpha_0,w}\\
 ~ & = & \lambda_1 S(g) D(\alpha ) D(\alpha_0) S(w) e_0\\
~ & = & \lambda_2 S(g) D(\alpha_1 ) S(w) e_0,
\end{eqnarray*}
where $\alpha_1 =\alpha +\alpha_0$ and 
$\lambda_2=\lambda_1e^{i\theta_h(\alpha_1,\alpha_0)}$. With
equations  (\ref{maimult3}),
(\ref{alpha2}), we have $I=\lambda_2 D(\alpha_2)S(g)S(w)e_0$, where
$\alpha_2=a\alpha_1+b\bar{\alpha}_1$. But  (\ref{r1}) implies
$I=\lambda_3 D(\alpha_2)S(g)e_{0,w}$, with
$\lambda_3=\lambda_2(1-|w|^2)^k$. Now we use  (\ref{r3}) and we
find $I=\lambda_4 D(\alpha_2)e_{0,w_1}$, where in accord with
 (\ref{r4}) $w_1$ is given by
(\ref{xxx}) and $\lambda_4=(\bar{a}+\bar{b}w)^{-2k}\lambda_3$. We
rewrite the last equation as $I= \lambda_5 D(\alpha_2)S(w_1)e_0$,
where $\lambda_5= (1-|w_1|^2)^{-k}\lambda_4$. Then we apply again
Lemma \ref{lema6} and we find $I=\lambda_6e_{z_1,w_1}$, where
$\lambda_6=\lambda_5(1-|w_1|^2)^k\exp (-\frac{\bar{\alpha}_2}{2}{z_1})$,
and  $z_1=\alpha_2-w_1\bar{\alpha}_2$. Proposition \ref{mm1} is
proved. \hfill  $\gata$

\begin{Remark}\label{rem9} Combining the expressions
 {\em{(\ref{xxx})-(\ref{x3})}},   the
factor $\lambda$ in  {\em{(\ref{xx})}} can be written down as
\begin{equation}\label{x4}
\lambda= (\bar{a}+\bar{b}w)^{-2k}
\exp (-\lambda_1), 
\end{equation}
where 
\begin{equation}\label{x5}
\lambda_1=\frac{\bar{b}z^2+
(\bar{a}\bar{\alpha}+\bar{b}\alpha)(2z+z_0)}{2(\bar{a}+\bar{b}w)}, 
~ z_0= \alpha-\bar{\alpha} w,
\end{equation}
or
\begin{equation}\label{x6}
\lambda_1=\frac{\bar{b}(z+z_0)^2}{2(\bar{a}+\bar{b}w)}
+\bar{\alpha}(z+\frac{z_0}{2}).  
\end{equation}
Note the expression {\em{(\ref{x4})-(\ref{x6})}} is identical with the
expression given in {\em Theorem  1.4} in {\em \cite{ez}}
 of the Jacobi forms, under
the the identification of $c, d,\tau, z, \mu, \lambda$ in \cite{ez}
with, respectively,  $\bar{b}, \bar{a}, w, z,\alpha ,-\bar{\alpha}$ in our
notation. Note also that the composition law {\em{(\ref{compositie})}} of the
Jacobi group $G^J$ and the action of the Jacobi group on the base
manifold {\em{(\ref{mm})}} is similar with that 
 in the paper \cite{bb}. {\em See also  \S \ref{rem19} 
and the Corollary 3.4.4 in \cite{bs}.}
\end{Remark}
\section{The symmetric Fock space}\label{sFs}

We recall the construction (\ref{aa}) of the map
$$\Phi :\Hi^*\rightarrow \fl; \Phi (\psi )=f_{\psi},
~f_{\psi}(z):=(e_{\bar{z}},\psi )_{\Hi}, $$
and the isometric embedding (\ref{anti}). 

Knowing  the symmetric
Fock spaces associated to the groups HW and $\text{SU}(1,1)$,
 we shall construct in this section  the symmetric Fock space associated to
the Jacobi group.

We begin recalling  the construction  for
\subsection{The Heisenberg-Weyl group}\label{heil}

In the orthonormal base (\ref{ort2}),
  Perelomov's CS-vectors associated to the HW
group, defined on $M:=H_1/\R=\C$, are
\begin{equation}\label{v1}
 e_z:=e^{za^+}e_0=\sum\frac{z^n}{(n!)^{1/2}}|n> ,
\end{equation}
and their corresponding holomorphic functions  are (see e.g. \cite{bar})
\begin{equation}\label{f1}
f_{|n>}(z):=(e_{\bar{z}},|n>)=\frac{z^n}{(n!)^{1/2}} .
\end{equation} 

The reproducing kernel $K:\C\times\bar{\C}\rightarrow\C$ is
\begin{equation}\label{ker1}
K(z,\bar{z}'):=(e_{\bar{z}},{e}_{\bar{z}'})=\sum f_{|n>}(z)\bar{f}_{|n>}(z')=
e^{z\bar{z}'} ,
\end{equation}
where the vector $e_z$ is given by  (\ref{v1}), while the function
$f_{|n>}(z)$ is given by  (\ref{f1}). In order to obtain the equality (\ref{ker1})
with $e_z$ given by  (\ref{v1}),  equation  (\ref{ort2}) is used. 

The scalar product (\ref{scf})
on the Segal-Bargmann-Fock space  is (cf. \cite{bar}) 
$$(\phi ,\psi )_{\Hi^*}=(f_{\phi},f_{\psi})_{\fl}=\frac{1}{\pi}
\int \bar{f}_{\phi}(z)f_{\psi}(z)e^{-|z|^2}d\Re z d\Im z .$$

Now we recall the  similar construction for 
\subsection{The group $\text{SU}(1,1)$}\label{su1}
In the  orthonormal base  (\ref{ort1}), 
Perelomov's CS-vectors for $\text{SU}(1,1)$, based on the unit disk
$\mc{D}_1=\text{SU}(1,1)/\text{U}(1)$,  are
\begin{equation}\label{v2}
e_z:=e^{z{\mb{K}}_+}e_0= \sum \frac{z^n{\mb{K}}_+^n}{n!}e_0= 
\sum \frac{z^n e_{k,k+n}}{n!a_{kn}} ,\end{equation}
and  the corresponding holomorphic functions 
are (see e.g.   equation 9.14 in \cite{bar47})
\begin{equation}\label{f2}
f_{e_{k,k+n}}(z):=(e_{\bar{z}},e_{k,k+n})
=\sqrt{\frac{\Gamma (n+2k)}{n!\Gamma (2k)}}z^n .
\end{equation}

The reproducing kernel $K:\mc{D}_1\times\bar{\mc{D}}_1\rightarrow \C$ is
\begin{equation}\label{ker2}
K(z,\bar{z}'):=(e_{\bar{z}},{e}_{\bar{z}'})=
\sum f_{e_{k,k+m}}(z)\bar{f}_{e_{k,k+m}}(z')=
(1-z\bar{z}')^{-2k}, 
\end{equation}
where the vector $e_z$ is given by  (\ref{v2}), while the function
$f_{e_{k,k+m}}(z)$ is given by  (\ref{f2}). In order to obtain the 
 equality  (\ref{ker2}) for $e_z$ given by  (\ref{v2}),
 the orthonormality  given by 
 (\ref{ort1}) is used, while for the second equality
involving the functions (\ref{f2}), use is made of
equation  (\ref{rajik}).  

The scalar product (\ref{scf}) on
$\mc{D}_1=\text{SU}(1,1)/\text{U}(1)$ 
 is  (see e.g.
equation 9.9 in \cite{bar47})
$$(\phi ,\psi )_{\Hi^*}=(f_{\phi},f_{\psi})_{\fl}=\frac{2k-1}{\pi}
\int_{|z|<1} \bar{f}_{\phi}(z)f_{\psi}(z)(1-|z|^2)^{2k-2}d\Re zd\Im z
, 2k= 2,3,... .$$

\subsection{The Jacobi group}\label{jcg}

In  formula (\ref{csu}) defining  Perelomov's CS vectors for the
 Jacobi group  (\ref{jac}), we take into account  (\ref{kkk1}),
(\ref{kkk3}) and we have
$$e_{z,w}=\exp (za^+ +\frac{1}{2}(a+)^2w)\exp (w{\mb{K}}'_+)e_0 .$$
With  (\ref{cucu1}), (\ref{ee}),  we have
$$e_{z,w}=\sum_{n}\frac{i^{-n}}{n!}
(\frac{w}{2})^{\frac{n}{2}}H_n(\frac{iz}{\sqrt{2w}})
(a^+)^n\sum_{m}\frac{w^m}{m!}({\mb{K}}'_+)^me_0 .$$
Now we take into account  (\ref{ort1}) and we get
 $$e_{z,w}=\sum_{n}\frac{i^{-n}}{(n!)^{1/2}}|n>
(\frac{w}{2})^{\frac{n}{2}}H_n(\frac{iz}{\sqrt{2w}})
\sum_{m}\frac{w^m}{m!a_{k'm}}e_{k',k'+m} .$$
The base of functions associated to the CS-vectors attached  to the
Jacobi group (\ref{jac}), based on the manifold $M$ (\ref{mm}) 
\begin{equation}\label{x3x1}
f_{{|n>;e_{k',k'+m}}}(z,w)  := 
(e_{\bar{z},\bar{w}},|n> e_{k',k'+m}),~z\in\C,~ |w|<1, 
\end{equation}
consists 
of the functions
\begin{equation}\label{x3x}
 f_{{|n>;e_{k',k'+m}}}(z,w) = 
(n!)^{-1/2}(\frac{i}{\sqrt{2}})^n\sqrt{\frac{\Gamma
(m+2k')}{m!\Gamma (2k')}}w^{m+\frac{n}{2}}H_n(\frac{-iz}{\sqrt{2w}}).
\end{equation}
Using the equation (\ref{hermite}), we can write down
\begin{equation}\label{x6x}
(\frac{i}{\sqrt{2}})^n
w^{\frac{n}{2}}H_n(\frac{-iz}{\sqrt{2w}}):=P_n(z,w) ,
\end{equation} where the polynomials $P_n(z,w)$ have the expression
\begin{equation}\label{marea}
P_n(z,w)=n!\sum _{k=0}^{[\frac{n}{2}]}
(\frac{w}{2})^k\frac{z^{n-2k}}{k!(n-2k)!} .
\end{equation} 
With the notation (\ref{x6x}), equation (\ref{x3x}) can be written down as
\begin{equation}\label{x4x}
f_{{|n>;e_{k',k'+m}}}(z,w)=f_{e_{k',k'+m}}(w)\frac{P_n(z,w)}{\sqrt{n!}}, 
\end{equation}
where the functions $f_{e_{k',k'+m}}$ are defined in  (\ref{f2}).

Above we have $2k=2k'+1/2$ and $m=0,1,\ldots $.

In order to illustrate   (\ref{marea}), we present the first 6
polynomials $P_n(z,w)$:
\begin{equation}
\begin{array}{ll}
P_0= 1; & P_1 = z ;\\
P_2 = z^2 +w ; & P_3 = z^3 + 3 zw ; \\
P_4= z^4 + 6 z^2w +3 w^2 ; & P_5 = z^5 +10 z^3 w + 15 z w^2 .
\end{array}
\end{equation}
 
The reproducing kernel (\ref{KHK}) $K:M\times \bar{M}\rightarrow\C$
 has the property:
\begin{subequations}\label{ker3}
\begin{eqnarray}
\label{sub1}K(z,w;\bar{z},\bar{w}')
 & := & (e_{\bar{z},\bar{w}},e_{\bar{z}',\bar{w}'})
= \sum_{n,m}f_{|n>,e_{k',k'+m}}(z,w)\bar{f}_{|n>,e_{k',k'+m}}(z',w') \\
\label{sub2} & = & (1-{w}\bar{w}')^{-2k}\exp{\frac{2\bar{z}'{z}+z^2\bar{w}'
+\bar{z}'^2w}{2(1-{w}\bar{w}')}} .
\end{eqnarray}
\end{subequations}
The fact that the coherent state vectors (\ref{csu}) have the scalar
product given by  (\ref{sub2}) was proved in Lemma 
\ref{lema5}, equation (\ref{KHK}). In order to check up the equality
(\ref{ker3})  for
the functions (\ref{x3x}), we use the form (\ref{x4x}), sum up the
part corresponding to the functions (\ref{f2}) using the summation
formula
(\ref{ker2}) for $k'$, while for the part corresponding to the mixed
part $P_n(z,w)$ in  (\ref{x4x}), we go back  with  (\ref{x6x}) to the
representation (\ref{x3x}), and we apply again the summation formula 
(\ref{bat})
of the Hermite polynomials.

In accord with the general scheme of \S \ref{scf}, the scalar product
(\ref{scf}) of
functions from the space $\mc{F}_K$ corresponding to the kernel
defined by (\ref{KHK})  on the manifold (\ref{mm}) is:
\begin{equation}\label{ofi}
(\phi ,\psi )= \Lambda\! \int_{z\in\C;
|w|<1}\!\bar{f}_{\phi}(z,w)f_{\psi}(z,w)
(1\!-\!w\bar{w})^{2k}\!\exp{\!-\frac{|z|^2}{1\!-\!w\bar{w}}}
\exp{\!-\frac{z^2\bar{w}\!+\!\bar{z}^2w}{2(1\!-\!w\bar{w})}} d \nu ,
\end{equation}
where the value of the $G^J$-invariant measure $d\nu $ 
\begin{equation}\label{ofi3}
d \nu =\frac{d \Re w d\Im w}{(1-w\bar{w})^3}d \Re z d \Im z
\end{equation}
will be deduced latter in  (\ref{dnu}), in accord with the receipt
  (\ref{scf1}).

In order to find the value of the constant $\Lambda$ in
 (\ref{ofi}), we take the functions $\phi,\psi=1$, we change the variable
$z\rightarrow (1-w\bar{w})^{1/2}z$ and we get
$$1 = \Lambda \int_{|w|<1}(1-w\bar{w})^{2k-2}d \Re w d\Im w
\int_{z\in\C}\exp({-|z|^2})\exp(-\frac{z^2\bar{w}+\bar{z}^2w}{2})
d \Re z d \Im z. $$ 

We apply equations (A1), (A2) in \cite{bar70}:
$$I(B,C)=\int \exp
(\frac{1}{2}(z.Bz+\bar{z}.C\bar{z}))\pi^{-n}e^{-|z|^2}
\prod_{k=1}^n d\Re z_k d\Im z_k= [\det (1-CB)]^{-\frac{1}{2}}, $$
where $B$, $C$ are complex symmetric matrices such that $|B|<1,
|C|<1$. Here $n=1$, $B=-\bar{w}$, $C= -w$. So, we get
$$1= \pi \Lambda \int_{|w|<1}(1-\bar{w}w)^{2k-5/2}d \Re w d \Im w .$$
But $$\int_{|w|<1}\frac{d \Re w d\Im
w}{(1-|w|^2)^{\lambda}}=\frac{\pi}{1-\lambda },~{\text{where}}~ \lambda <1 ,$$
and we find out the value of the constant $\Lambda$ in  (\ref{ofi}):
\begin{equation}\label{ofi1}
\Lambda = \frac{4k-3}{2\pi^2} .
\end{equation}

\subsection{The geometry of the manifold $\C\times\mc{D}_1$}\label{two}
Now we follow the general prescription of \S\ref{s22}.
We calculate the K\"ahler potential
as the logarithm of the reproducing kernel ({\ref{KHK}}),
$f:=\log K$, 
i.e.
\begin{equation}\label{keler}
f =\frac{2z\bar{z}+z^2\bar{w}+\bar{z}^2w}{2(1-w\bar{w})}
-2k\log (1-w\bar{w}) .
\end{equation}
The K\"ahler two-form $\omega$ is given by
the formula:
\begin{equation}\label{aaa}
-i \omega = f_{z\bar{z}}dz\wedge d\bar{z}+f_{z\bar{w}}dz\wedge
d\bar{w}
-f_{\bar{z}w}d\bar{z}\wedge d w +f_{w\bar{w}}dw\wedge d\bar{w} .
\end{equation}
The volume form is:
\begin{equation}\label{vol1}-\omega\wedge\omega=2\left|\begin{array}{cc}
f_{z\bar{z}} & f_{z\bar{w}}\\
f_{\bar{z}w} & f_{\bar{w}w}
\end{array}\right|
dz\wedge d \bar{z}\wedge d w \wedge d \bar{w} .
\end{equation}
We start  calculating the partial derivatives of the function $f$. We have
$$f_{\bar{z}}=\frac{z+\bar{z}w}{1-w\bar{w}},$$
$$f_{\bar{w}}=
\frac{z^2+w(2z\bar{z}+\bar{z}^2w)}{2(1-w\bar{w})^2}+
\frac{2k w}{1-w\bar{w}},$$
$$f_{w\bar{w}}=\frac{A''}{2(1-w\bar{w})^4}+\frac{2k}{(1-w\bar{w})^2},$$
where $A''$ is
$$A''=(2\bar{z}^2w+2z\bar{z})(1-w\bar{w})^2+
2\bar{w}(1-w\bar{w})(\bar{z}^2w^2+2z\bar{z}w+z^2)=(1-w\bar{w})A' .$$
Here $A'$ is
$$A'=2\bar{z}^2w+2z\bar{z}-2\bar{z}^2w^2\bar{w}-2z\bar{z}\bar{w}w+
2\bar{z}^2\bar{w}w^2+4z\bar{z}w\bar{w}+2z^2\bar{w} .$$
and we have finally
$$A'=2(\bar{z}+\bar{w}z)(z+w\bar{z}) .$$

So, {\it we find for the manifold} (\ref{nmm})
  {\it the fundamental two-form $\omega$} (\ref{aaa}), {\it where}
\begin{subequations}\label{aas}
\begin{eqnarray}
f_{z\bar{z}} & = & \frac{1}{1-w\bar{w}}\label{aa1}, \\
f_{z\bar{w}}& = & \frac{z+w\bar{z}}{(1-w\bar{w})^2}\label{aa2} ,\\
f_{w\bar{w}} & = &
\frac{(\bar{z}+\bar{w}z)(z+w\bar{z})}{(1-w\bar{w})^3}
+\frac{2k}{(1-w\bar{w})^2} .\label{aa3}
\end{eqnarray}
\end{subequations}
We can write down the two-form $\omega$ (\ref{aaa})-(\ref{aas}) as
\begin{equation}\label{aab}
-i\omega =\frac{2k}{(1-w\bar{w})^2}dw \wedge d\bar{w} +
\frac{A\wedge \bar{A}}{1-w\bar{w}},
~A=dz+\bar{\alpha}_0dw, ~\alpha_0=\frac{z+\bar{z}w}{1-w\bar{w}}.
\end{equation}
For the volume form (\ref{vol1}), we find:
\begin{equation}\label{dnu}
\omega\wedge\omega = 4k (1-w\bar{w})^{-3} 4\Re z\Im z \Re w \Im w .
\end{equation}

It can be checked up that indeed, {\it  the measure
$d\nu$ and the fundamental two-form $\omega$ are group-invariant at the action 
{\em (\ref{xxx})} of the Jacobi group {\em (\ref{jac})}}.

Now we summarize the contents of this section as follows:
\begin{Proposition}\label{final}
Let us consider the Jacobi group $G^J_1$ {\em (\ref{jac})}  with the composition
rule
{\em (\ref{compositie})} acting on the
coherent state manifold {\em (\ref{nmm})} via equation 
 {\em  (\ref{xxx})}. The manifold
$\mc{D}^J_1$ has the K\"ahler potential {\em (\ref{keler})}  and the
$G^J_1$-invariant K\"ahler two-form $\omega$ given by {\em (\ref{aab})}.
 The holomorphic polynomials {\em(\ref{x3x1})} associated to the coherent
state vectors {\em (\ref{csu})} are given by {\em  (\ref{x4x})}, where the
functions
 $f$ are
given by  {\em (\ref{f2})}, while the polynomials $P$ are given by
{\em  (\ref{marea})}. The Hilbert space of holomorphic functions $\mc{F}_K$
associated to the holomorphic kernel $K:\mc{D}^J_1
\times \bar{\mc{D}}^J_1\rightarrow\C$ given by 
 {\em (\ref{KHK})} is endowed with the
scalar product {\em (\ref{ofi})}, where the normalization
 constant $\Lambda$ is given by
{\em (\ref{ofi1})} and the $G^J_1$-invariant measure $d\nu$
 is given by {\em (\ref{ofi3})}.
\end{Proposition}

Recalling Proposition IV.1.9. p. 104  and Proposition
XII.2.1 p. 515 in \cite{neeb}, Proposition
\ref{mm1} can be formulated as follows:
\begin{Proposition}\label{finalf}
Let $h := (g,\alpha )\in G^J_1$,  where $G^J_1$ is the Jacobi group
{\em (\ref{jac})}, and we consider the representation
  $\pi (h) := S(g)D(\alpha )$, $g\in
\rm{SU}(1,1)$,
$\alpha\in \C$, and let the notation $x:=(z,w)\in
{\mc{D}^J_1}:=\C\times\mc{D}_1$. Then the 
continuous unitary representation $(\pi_K,\Hi_K)$  attached to the
positive definite holomorphic kernel $K$ defined by {\em  (\ref{KHK})} is
\begin{equation}\label{rep}
(\pi_K(h).f)(x)=J(h^{-1},x)^{-1}f(h^{-1}.x),
\end{equation}
where the cocycle $J(h^{-1},x)^{-1}:=\lambda (h^{-1},x)$  with
$\lambda$ defined by equations
{\em  (\ref{xx})-(\ref{x3})} and the function $f$
belongs to the Hilbert space of holomorphic functions $\Hi_K\equiv
\mc{F}_K$ endowed with the
scalar product {\em (\ref{ofi})}, where $\Lambda$ is given by 
{\em  (\ref{ofi1})}.
\end{Proposition}

\begin{Remark}
The equations of the geodesics on the manifold $\mc{D}^J_1$
endowed with the two-form {\em({\ref{aab}})} in the variables $w\in\mc{D}_1, z\in\C$
are
\begin{subequations}
\begin{eqnarray}
2k\frac{d^2z}{dt^2}-\bar{\alpha}_0(\frac{dz}{dt})^2+2(2k\frac{\bar{w}}{P}-
\bar{\alpha}^2_0)\frac{dz}{dt}\frac{dw}{dt}-\bar{\alpha}^3_0(\frac{dw}{dt})^2
&= & 0;\\
2k\frac{d^2w}{dt^2}+(\frac{dz}{dt})^2+2\bar{\alpha}_0\frac{dz}{dt}\frac{dw}{dt}
+ (4k\frac{\bar{w}}{P}+\bar{\alpha}^2_0)(\frac{dw}{dt})^2 & = & 0, 
\end{eqnarray}
\end{subequations}
where $\alpha_0$ is given by  {\em{(\ref{x2})}} and
$P=1-w\bar{w}$.\end{Remark}

\section{Physical applications: classical and quantum
equations of motion}\label{app}

We consider  applications of the formulas (\ref{summa}) proved in this paper
for the Jacobi group $G^J_1$ to equations of motion on the
CS-manifold $\mc{D}^J_1$.
 This extend our previous results for hermitian groups 
\cite{sbcag,sbl} or semisimple groups which generate CS-orbits
\cite{sbctim,sin}
 to an example of  non-reductive CS-group.

Passing on from the dynamical system problem
 in the Hilbert space $\Hi$ to the corresponding one on $M$ is called
sometimes {\it dequantization}, and the system on $M$ is a classical
one \cite{sbcag},\cite{sbl}. Following Berezin \cite{berezin1}, the
motion on the classical phase space can be described by the local
equations of motion
\begin{equation}\label{ecmiscare}
\dot{z}_{\gamma}=i\{\mb{\mc{H}} ,z_{\gamma}\},
\end{equation}
where $\mb{\mc{H}}$ is the energy function attached to the Hamiltonian
$\mb{H}$.
In  (\ref{ecmiscare}) $\{\cdot,\cdot \}$ denotes  the  Poisson bracket:
$$\{ f,g\} =\sum_{\alpha ,\beta\in\Delta_+}\NC^{-1}_{\alpha
,\beta}\left\{\frac{\pa f}{\pa
z_{\alpha}}\frac{\pa g}{\pa\bar{z}_{\beta}}-\frac{\pa f}{\pa\bar{z}_{\alpha}}
\frac{\pa g}{\pa z_{\beta}}\right\}, f,g\in C^{\infty}(M) ,$$
where $\NC$ is defined in  (\ref{twoform}).
The equations of motion (\ref{ecmiscare}) can be written down as
\begin{equation}
i\left(
\begin{array}{cc}
0 & \NC\\
 & \\
-\bar{\NC} & 0
\end{array}\right)
\left(
\begin{array}{c}
\dot{z} \\
  \\
\dot{\bar{z}}
\end{array}
\right)
=-\left(
\begin{array}{c}
\frac{\pa}{\pa z}\\
   \\
\frac{\pa}{\pa \bar{z}}
\end{array}
\right)\mb{\mc{H}} .
\end{equation}

We consider an algebraic Hamiltonian linear in the generators of the
group of symmetry 
\begin{equation}\label{hsum}
\mb{H}=\sum_{\lambda\in\Delta}\epsilon_{\lambda}{\mb{X}}_{\lambda} .
\end{equation}
We recall (cf. \cite{sbcag,sbl}) that {\it if the differential action
of the generators of the group $G$
is given by formulas {\em(\ref{sss})}, then 
  classical motion and the quantum evolutions generated
by the Hamiltonian {\em (\ref{hsum})} are given by the same equations
of motion {\em (\ref{ecmiscare})} on $M=G/H$}: 
\begin{equation}\label{move}
{i\dot{z}_{\alpha}=\sum_{\lambda}\epsilon_{\lambda}Q_{\lambda
,\alpha}}.
\end{equation}
 {\it The two-form form $\omega$ on $M$ permits to determine
the Berry phase} \cite{sbl}.

Let us consider a linear Hamiltonian in the generators of the Jacobi
group (\ref{jac}):
\begin{equation}\label{guru}
\mb{H} = \epsilon_aa +\bar{\epsilon}_aa^+
 +\epsilon_0 {\mb{K}}_0 +\epsilon_+{\mb{K}}_++\epsilon_-{\mb{K}}_-  .
\end{equation}
We have
\begin{Remark}
The equations of motion on the manifold {\em(\ref{mm})}
 generated by
the linear Hamiltonian {\em(\ref{guru})}
 are given by the matrix Riccati equation:
\begin{subequations}
\begin{eqnarray}
i\dot{z} & = & \epsilon_a+\frac{\epsilon_0}{2}z +{\bar{\epsilon}}_aw+\epsilon_+ z w ,\\
i\dot{w} & = & \epsilon_- +\epsilon_0w+
\epsilon_+w^2 .\label{guru1}
\end{eqnarray}
\end{subequations}
\end{Remark}
Note that the second equation (\ref{guru1})
 is a Riccati equation on $\mc{D}_1$. The procedure of linearization
of matrix Riccati equation on manifolds is discussed in \cite{sbl}.

An interesting development of the present construction should  be to consider
 nonlinear CSs \cite{siv} attached to a deformed Jacobi
group. However, difficulties in the physical interpretation of the
creation and annihilation of $q$-deformed oscillator, related to the
quantum groups $\text{SU}_q(2)$ \cite{bied}, appear as these are
symmetric but not self-adjoint operators \cite{kli}.

\section{Comparison with  K\"ahler-Berndt's approach}\label{rem19}

  Rolf Berndt -alone or in collaboration - has studied the 
 real Jacobi group $G^J(\R )$ in
several references, from which I mention 
\cite{bern,bb,bs,berndt}. The Jacobi group appears (see explanation in
\cite{cal}) in the
context of the so called {\it Poincar\'e group} or {\it The New
Poincar\'e group} investigated by Erich K\"ahler as the 10-dimensional
group $G^K$ which invariates a hyperbolic metric
\cite{cal1,cal2,cal3}. K\"ahler and Berndt have investigated the
Jacobi group  $G^J_0(\R ):= \text{SL}_2(\R )\ltimes \R^2$ acting on the
manifold $\mc{X}^J_1:=\mc{H}_1\times\C$, where  $\mc{H}_1$ is the upper half plane
$\mc{H}_1:=\{v\in\C|\Im (v)>0\}$. 

For self-contentedness, in
 Remarks \ref{actnea}  and \ref{actnea1} below, we briefly proof two
results from \cite{bs}, which we need in order two express our two
form $\omega$ in the coordinates used by K\"ahler and Berndt. 
\begin{Remark}\label{actnea}
The action of $G^J_0(\R )$ on $\mc{X}^J_1$ is
given by $(g,(v,z))\rightarrow (v_1,z_1)$, $g=(M , l )$, where
\begin{equation}\label{lac}
v_1=\frac{av+b}{cv+d},  z_1=\frac{z+l_1v+l_2}{cv+d}
; ~M  = \left(\begin{array}{cc} a & b\\c & d
\end{array}\right)\in \rm{SL}_2(\R) , (l_1,l_2)\in\R^2 .
\end{equation} 
\end{Remark}
{\it {Proof}}. 
Let us use the notation of \cite{bs}. We denote 
$G^J(\R ):=\SL\ltimes H(\R)$, where  $H(\R)$ is  the real HW group
with the composition law:
\begin{equation}\label{cplh}
(\lambda ,\mu, \ka )(\lambda',\mu',\ka')=(\lambda +\lambda',\mu
+\mu',\ka +\ka' +\left|\begin{array}{c} X\\X'\end{array}\right|)
, \left|\begin{array}{c} X\\X'\end{array}\right|= \det 
\left(\begin{array}{c} X\\X'\end{array}\right) .
\end{equation} 
  If $g=(M,X,\ka )\in G^J (\R )$, where $M\in \SL$, $X=(\lambda,\mu )$, $(X,\ka
)\in \R^3$, then the composition law in the
real Jacobi group  is
\begin{equation}\label{clj}
gg'=(MM',XM'+X', \ka +\ka '+\left|\begin{array}{c} XM'\\X'\end{array}\right|).
\end{equation}
The action of $G^J(\R )$ over the $H(\R)$ is
\begin{equation}\label{onheis}
M(X,\ka )M^{-1}=(XM^{-1},\ka ). 
\end{equation}
Let us consider the Iwasawa decomposition for  a matrix $M\in\SL$:
\begin{equation}\label{iwa}
M=\left(\begin{array}{cc}1 &x \\ 0 & 1 \end{array}\right)
\left(\begin{array}{cc} y^{1/2} & 0 \\ 0 & y^{-1/2} \end{array}\right)
\left(\begin{array}{cc}\cos\theta &\sin\theta \\ -\sin\theta &
\cos\theta
 \end{array}\right), ~y>0.
\end{equation}
If 
\begin{equation}\label{abcx}
 M= \left(\begin{array}{cc}a & b \\ c & d\end{array} \right),
\end{equation}
then we find for $x,y, \theta$ in (\ref{iwa})
\begin{equation}\label{xyt}
x=\frac{ac+bd}{d^2+c^2};~ y=\frac{1}{d^2+c^2};~ \sin\theta
=-\frac{c}{\sqrt{c^2+d^2}};~
\cos\theta
=\frac{d}{\sqrt{c^2+d^2}}.
\end{equation}
For $g=(M,X,\ka )\in G^J(\R)$, the EZ-coordinates  (Eichler-Zagier,
cf. the definition at p. 12 and p. 51
in \cite{bs}) are
 $(x,y,\theta,\lambda,\mu,\ka )$. Let $\tau=x+
iy\in\mc{H}_1$, $z=\xi +i\eta =p\tau + q$, where
\begin{equation}\label{pqxy}
(p,q)=XM^{-1}=(\lambda d-\mu c, -\lambda b+ \mu a )
\end{equation}  
Using the multiplication law (\ref{clj}), the Iwasawa decomposition
(\ref{iwa}) and the equations (\ref{xyt}), (\ref{pqxy}), we find
the action of $ G^J(\R )$  on the base $\mc{X}^J_1$ 
 \begin{equation}
g(\tau ,z) =(\frac{a\tau+b}{c\tau+d},\frac{z+\lambda \tau+\mu}{c\tau +d}),
\end{equation}
and Remark \ref{actnea} is proved.
\hfill $\gata$

Let us now  recall that 
\begin{equation}\label{xcs}C^{-1}\text{SL}_2(\R )C=\text{SU}(1,1),~\text{where}~
C=\left(\begin{array}{cc} i& i \\
-1& 1\end{array}\right) .
\end{equation}
If $M\in \SL$ is the matrix (\ref{abcx}), then, under the transformation
(\ref{xcs})
\begin{equation}\label{mstar}
M^*=C^{-1}MC=\left(\begin{array}{cc}\alpha&\beta\\ \bar{\beta} &
\bar{\alpha}\end{array}\right),~
\alpha,\beta\in\C, |\alpha|^2-|\beta|^2=1,
\end{equation}
where 
\begin{equation}\label{mmstar}
2\alpha = a+d+i(b-c);~ 2\beta=
a-d-i(b+c).
\end{equation}

Now we pass to the complex group  $G^J_{\C}=C^{-1}G^J(\R )C$.
We  recall  that the Jacobi group
 $G^J_{\C }$  is a 
{\it group of Harish-Chandra type},  (cf. e.g.  p. 514 in \cite{neeb};
see the
 definition in  Ch. III \S 5 in \cite{satake} and
Ch. XII.1 in \cite{neeb}).  Moreover, it  is well known that 
{\it the Jacobi algebra} (\ref{baza}) {\it is a  CS-Lie algebra}
  (cf. e.g. Theorem 5.2  in \cite{lis}).
 The correspondence between our notation and that of
Berndt-Schmidt  at p. 12 in \cite{bs} is as follows: $a^+, a, K_+, K_-,
1, K_0$ corresponds,  
respectively to: $ Y_+, Y_-, X_+, -X_-, -Z_0, \frac{1}{2}Z$. See also
our Remark \ref{rem18} below.   

We see that under the  transformation (\ref{xcs}), $g =(M,X,\ka )\in
\SL\ltimes H(\R )$ is twisted to $g^*=(M^*,X^*,\ka )$, where  $M^*$ is
given
by (\ref{mstar}), while, due to action (\ref{onheis}), $X^*= XC=
(i\lambda -\mu, i\lambda + \mu)$. 

Also the map (\ref{xcs}) induces a transformation of the bounded
domain $\mc{D}_1$ into the upper half plane $\mc{H}_1$ and
\begin{equation}\label{taustar}
\tau\in\mc{H}_1\mapsto\tau^*=C^{-1}(\tau )=\frac{\tau -i}{\tau
+i}\in\mc{D}_1 .\end{equation}
 The action $C^{-1}G^J_0(\R )C$ 
 descends on the basis as the biholomorphic map: 
$\check{C}^{-1}: \mc{X}^J_1:= \mc{H}_1\times \C \rightarrow
\mc{D}^J_1:=\mc{D}_1\times \C$: $(\tau, z)\mapsto (\tau^*,z^*)$. Here
$\tau*$ is given by (\ref{taustar}), while $z^*=p^*\tau^*+q^*$. So,
$(p,q)= (\lambda ,\mu ) M^{-1}$, and $(p^*,q^*)= (\lambda^* ,\mu^* )
{M^*}^{-1}$. But
$M^*=C^{-1}MC$, and $(p^*,q^*)= (p,q)C= (-q+ip, q+ip)$, and we get
$z^*=\frac{2iz}{\tau + i}$.
Note that at  p. 53 in \cite{bs} 
the   factor $2i$ in this formula 
 is missing.

In a different notation,  we have shown   
  that
 \begin{Remark}\label{actnea1}
The action $C^{-1}G^J_0(\R )C$, 
 descends on the basis as the biholomorphic map: 
$\check{C}^{-1}: \mc{X}^J_1:= \mc{H}_1\times \C \rightarrow
\mc{D}^J_1:=\mc{D}_1\times \C$:
\begin{equation}\label{noiec}
w=\frac{v-i}{v+i};~ z=\frac{2iu}{v+i},
w\in\mc{D}_1,~v\in\mc{H}_1,z\in\C .
\end{equation}
\end{Remark}

The $G^J_0(\R )$-invariant closed two-form considered by  K\"ahler-Berndt
 is:
\begin{equation}\label{ura1}
\omega '= \alpha 
\frac{dv\wedge d
\bar{v}}{(v-\bar{v})^2}+\beta\frac{1}{v-\bar{v}}B\wedge
 \bar{B}, ~B=
du-\frac{u-\bar{u}}{v-\bar{v}}dv,
v,u\in\C, ~\Im (v)>0, 
\end{equation}
cf.\S 36 in \cite{cal3}; see also   \S 3.2 in \cite{bern}, where the
first term is misprinted as $\alpha 
\frac{dv\wedge d \bar{v}}{v-\bar{v}}$. 

 Under the mapping (\ref{noiec}), the two-form $\omega$ (\ref{aab}) reads
\begin{equation}\label{ura3}
-i\omega = -\frac{2k}{(\bar{v}-v)^2}dv\wedge d\bar{v}+\frac{2}{i(\bar{v}-v)}B\wedge \bar{B},
\end{equation}
i.e. (\ref{ura1}). In fact, we have proved that
\begin{Remark}
When expressed in the  coordinates $(v, u)\in\mc{X}^J_1$ which are related to the
coordinates $(w,z)\in\mc{D}^J_1$ by the map  
{\em{(\ref{noiec})}}  given by {\em Remark \ref{actnea1}},
our K\"ahler two-form {\em{(\ref{aab})}} is 
identical with the one  {\em({\ref{ura3}})}
 considered by K\"ahler-Berndt.
\end{Remark}

If we use the EZ coordinates 
adapted to our notation
\begin{equation}\label{eqez}
v=x+iy; ~ u=pv+q,~ x,p,q, y\in\R ,y>0,
\end{equation}
the $G^J_0(\R )$-invariant K\"ahler metric on $\mc{X}^J$
 corresponding to the K\"ahler-Berndt's
K\"ahler two-form $\omega$ (\ref{ura3}) reads
\begin{equation}\label{victorie}
ds^2=\frac{k}{2y^2}(dx^2+dy^2)+\frac{1}{y}[(x^2+y^2)dp^2+dq^2+2xdpdq],
\end{equation}
i.e. the equation at p. 62 in \cite{bs} or the equation
given  at p. 30 in \cite{bern}.

 The K\"ahler two-form (\ref{ura1}) of K\"ahler-Berndt corresponds
(cf. equation (9) in Ch. 36 of \cite{cal1})
to the K\"ahler potential
\begin{equation}\label{cpot}
f'=-\frac{\lambda}{2}\log
\frac{v-\bar{v}}{2i}-i\pi\mu\frac{(u-\bar{u})^2}{v-\bar{v}}, 
u\in\C, ~v\in\mc{H}_1.
\end{equation} 

The K\"ahler potential (\ref{cpot}), corresponds to some
{\it (K\"ahler) Perelomov's CS-vectors 
  based on the
CS-manifold $\mc{X}^J_1$} on which the group $G^J_0(\R)$ acts via the the action
(\ref{lac}), which, instead of the scalar product $K$ (\ref{hot}),
 should have {\it  a scalar product $K'$ in the EZ coordinates 
{\em({\ref{eqez}})}} $x,y,p,q\in\R$
\begin{equation}\label{hot11}
K'= y^{-\frac{\lambda}{2}} \exp ({2\pi \mu p^2y}).  
\end{equation}

It would be interesting to extend eq. (\ref{ura3}) to the case
${H}_n\rtimes\text{Sp}(n,\R )$, see reference \cite{sbj}. This would be
useful for a  better understanding of  the Gaussons \cite{si,ali}.

\section{Some more comments}\label{lastc}
\begin{Remark}\em
 We have called the algebra (\ref{baza}), the {\it Jacobi
algebra} and the
group (\ref{jac}), the {\it Jacobi group}, in agreement  with the name
used   in \cite{bs} or at p. 178 in \cite{neeb}, 
where the algebra $\got{g}^J_1:= 
\got{h}_1\rtimes\got{sl}(2,\R )$  is called ``Jacobi  algebra''.
 The denomination adopted in the present paper is of course
in accord with the one used in \cite{neeb}   because of the
isomorphism of the Lie algebras $\got{su}(1,1)
\sim\got{sl}(2,\R)\sim\got{sp}(1,\R )$$(\sim \got{so}(2,1))$. 
 Also the name ``Jacobi
algebra'' is used in \cite{neeb} p. 248 to call the semi-direct sum of
the $(2n+1)$-dimensional Heisenberg algebra and the symplectic algebra,
 $\got{hsp}:=
\got{h}_n\rtimes \got{sp}(n,\R )$. 
 The group corresponding to this algebra
 is called sometimes  in 
the  Mathematical Physics literature (see e. g. \S 10.1 in \cite{ali},
which is based on \cite{si})
the   ``metaplectic group'', but in reference \cite{neeb} the term 
``metaplectic group''  is reserved to  the 2-fold covering group of the
symplectic group, cf. p. 402 in  \cite{neeb} (see also \cite{bar70}
and \cite{itzik}). Other names of the metaplectic representation are
the oscillator representation, the harmonic representation or the
Segal-Shale-Weil representation, see references in Chapter 4 of
\cite{fol} and \cite{bs}. 
\end{Remark}

\begin{Remark}\label{rem18}\em
 We now discuss the Jacobi algebra (\ref{baza})
 from the view point of the book \cite{neeb}. 
We know (cf. Lemma XII.1.20 p. 509) that {\it quasihermitian}
 Lie algebras, i.e. Lie
algebras for which a {\it maximal compactly embedded subalgebra}
$\got{k}$ (cf. Definition VII.1.1
p. 222 in \cite{neeb}) satisfies the relation
$\got{z}_{\got{ g}}(\got{z}(\got{k}))=\got{k}$ (cf. Definition VII.2.15
p. 241 in \cite{neeb}), admit the 5-grading of the complexification
$\got{g}_{\C} =
\got{p}^+_s\oplus\got{p}^+_r\oplus\got{k}_{\C}\oplus\got{p}^-_r\oplus
\got{p}^-_s$, where $\got{k}$ is the maximal compactly embedded
subalgebra of $\got{g}$, $\got{p}_s$ ($\got{p}_r$) represent the
semisimple roots, (respectively, the solvable roots) (cf. Definition
VII.2.4 p. 234 in \cite{neeb}), while ``$^+$'' (``$^-$'') refers to the
positive (respectively, negative) roots
with respect to a $\Delta^+$ adopted
positive system (cf. Definition 
VII.2.6 p. 236 in \cite{neeb}). But {\it the Jacobi algebra is quasihermitian},
cf. Example VIII.2.3  p. 294 in \cite{neeb}, with 
$\got{t}=\{0\}\oplus\R\oplus\got{u}(1)$ a compactly embedded Cartan
algebra and also a maximal compactly  embedded subalgebra of $\got{g}^J$
(cf. Example VII.2.30 p. 249 and Example XII.1.22 p. 513   in
\cite{neeb}). So, the generators $K_+$, $K_-$, $a^+$, $a$ of the
Jacobi algebra (\ref{baza}) belong to
$\got{p}^+_s$, $\got{p}^-_s$, $\got{p}^+_r$, respectively $\got{p}^-_r
$,  $1$ belongs to the $\R$-part of $\got{k}$, while
$K_0$ belongs to the $\got{u}(1)$-part of $\got{k}$.
 Note that due to relation (\ref{baza4}), the
subalgebra $\got{p}^+=\got{p}^+_r\oplus\got{p}^+_s$ which appears in
the definition (\ref{csu}) of Perelomov's coherent state vectors is
an abelian one, as it should be
 (cf. Lemma XII.1.20 p. 509 in \cite{neeb}).
\end{Remark}
\begin{Remark}\em
 We emphasize that 
 {\it the representation given in  {\em Lemma \ref{mixt}}
is different from the extended metaplectic  representation}
(cf. p. 196 in
 \cite{fol}, see also \cite{satake}). As was already mentioned, the
Jacobi algebra admits a realization as subalgebra of the Weyl algebra $A_1$
of polynomials in $p,q$ of degree $\le 2$. In the present paper we
have presented a realization of the Jacobi algebra as subalgebra of
the Weyl algebra $A_2$ defined by holomorphic first order
differentials
 operators with holomorphic polynomials of degree $\le 2$. This
algebra is realized in the variables $x=(z,w)\in
\mc{D}^J_1$. We recall that the only finite-dimensional
non-solvable Lie algebras that can be realized as Lie subalgebras of
the complex Weyl algebra $A_1$ are: $\got{sp}(1,\R)$,
$\got{sp}(1,\R)\times \C$ and the Jacobi algebra \cite{sz}, \cite{jos},
\cite{rom}.
\end{Remark}
\begin{Remark}\em
  Note that  the expression   (\ref{KHK}) of
 the reproducing kernel is the  particular case $n=1$, $2k =\frac{1}{2}$
of  the reproducing kernel on the
space $\mc{D}^J_n:= \C^n\times \mc{D}_n$, where $\mc{D}_n$
is the Siegal ball, at  p. 532 in the book of
K.-H. Neeb \cite{neeb} or in the article \cite{hn} and
 (5.28)
 in the book
 \cite{satake}. See also \cite{sbj}. 
\end{Remark}

\section{Appendix}\label{APEN}

For self-completeness, we also give
\vspace*{1cm}

{\it {Proof  of {\em  (\ref{ort2})}}}.
We shall calculate 
$$\lambda_{n;m}:=(e_0,a^n(a^+)^me_0).$$
We shall use the formula
\begin{equation}\label{cheia}
[A,B^m]=\sum_{s=0}^{m-1}B^s[A,B]B^{m-s-1},
\end{equation}
for $A=a^n$, $B=a^+$. We have
$$[A,B^m]=-\sum_{s=0}^{m-1}B^s[B,A]B^{m-s-1}.$$
But
$$[B,A]=[a^+,a^n]=\sum_{p=0}^{n-1}a^p[a^+,a]a^{n-p-1}=-\sum_{p=0}^{n-1}a^{n-1}=
-na^{n-1},$$
and
$$[A,B^m]= n\sum_{s=0}^{m-1}(a^+)^sa^{n-1}(a^+)^{m-s-1},$$
so, we have
$$\lambda_{n;m}=n\lambda_{n-1;m-1}.$$
If $n=m$, then $\lambda_{nn}=n!$.
If $n<m$, then $\lambda_{n;m} = ct(e_0,[a,(a^+)^p]e_0)$, where $p>1$.
But
$[a,(a^+)^p]=\sum_{s=0}^{p-1}(a^+)^s[a,a^+](a^+)^{p-s-1}=p(a^+)^{p-1}$,
and
$\lambda_{n;m}= ct (e_0,(a^+)^{p-1}e_0)=0$.

Similarly, if $n>m$, then $[a^p,a^+]=pa^{p-1}$ and also $\lambda_{n;m}=
ct (e_0,a^{p-1}e_0)=0$.

So, we have $\lambda_{n;m}=n!\delta_{n;m}$.  \hfill $\gata$
\vspace*{1cm}

{\it {Proof  of {\em  (\ref{ort1})}}}.
We calculate $\mu_{n;m}:=(e_0,C_{n;m}e_0)$, where
$C_{n;m}={\mb{K}}_-^n{\mb{K}}_+^m$, using  (\ref{cheia}) with
 $A={\mb{K}}^n_-$   and  $B={\mb{K}}_+$. We find
$$\mu_{n;m}= (e_0,\sum_{s=0}^{m-1}{\mb{K}}^s_+
[{\mb{K}}^n_-,{\mb{K}}_+]{\mb{K}}^{m-s-1}_+e_0)= 
(e_0,[{\mb{K}}^n_-,{\mb{K}}_+]{\mb{K}}^{m-1}_+e_0) .$$
But
\begin{eqnarray*}
[{\mb{K}}_+,{\mb{K}}_-^n] & =&
\sum_{p=0}^{n-1}{\mb{K}}^p_-[{\mb{K}}_+,{\mb{K}}_-]{\mb{K}}^{n-p-1}_-
=  
-2\sum_{p=0}^{n-1}{\mb{K}}^p_-{\mb{K}}_0{\mb{K}}^{n-p-1}_-\\
  & = & -2\sum_{p=0}^{n-1}{\mb{K}}^p_-[{\mb{K}}_0,{\mb{K}}^{n-p-1}_- ]
-2\sum_{p=0}^{n-1}{\mb{K}}^p_-{\mb{K}}^{n-p-1}_-{\mb{K}}_0\\
& = & -2n{\mb{K}}^{n-1}_-{\mb{K}}_0 -
2\sum_{p=0}^{n-1}{\mb{K}}^p_-[{\mb{K}}_0,{\mb{K}}_-^{n-p-1}].
\end{eqnarray*}
We find
$$\mu_{n;m}= 2n (e_0,{\mb{K}}_-^{n-1}{\mb{K}}_0{\mb{K}}^{m-1}_+e_0)+R,
$$
where
$$R:=2(e_0,R_0e_0); ~
R_0:=
\sum_{p=0}^{n-1}{\mb{K}}^p_-[{\mb{K}}_0,{\mb{K}}_-^{n-p-1}]{\mb{K}}^{m-1}_+
,$$
and we get 
$$\mu_{n;m}=2nk\mu_{n-1,m-1}+
2n(e_0,{\mb{K}}_-^{n-1}[{\mb{K}}_0,{\mb{K}}^{m-1}_+]e_0)+R .$$
But 
$$[{\mb{K}}_0,{\mb{K}}^{m-1}_+]=\sum_{k=0}^{m-2}{\mb{K}}^{s}_+
[{\mb{K}}_0,{\mb{K}}_+]{\mb{K}}^{m-2-s}_+=(m-1){\mb{K}}^{m-1}_+,$$
and
$$[{\mb{K}}_0,{\mb{K}}^{n-p-1}_-]=
\sum_{q=0}^{n-p-2}{\mb{K}}^q_-[{\mb{K}}_0,{\mb{K}}_-]{\mb{K}}_-^{n-p-q-2}=
-(n-p-1){\mb{K}}^{n-p-1}_-.$$
We get successively
$$R_0=-\sum_{p=0}^{n-1}(n-p-1)C_{n-1;m-1},$$
$$R=-n(n-1)\mu_{n-1;m-1},$$
and
$$\mu_{n;m}=(2nk + 2n(m-1) -n(n-1))\mu_{n-1;m-1},$$
$$\mu_{n;n}=n(2k+n-1)\mu_{n-1;n-1}; ~\mu_{1;1}=2k,$$
$$\mu_{n;n}=\frac{n!(2k+n-1)!}{(2k-1)!}=\frac{n!\Gamma
(2k+n)}{\Gamma (2k)}.$$

If $n<m$, then there is a $p>1$ such that
$$[{\mb{K}}_-,{\mb{K}}_+^p]=\sum_{q=0}^{p-1}{\mb{K}}^q_+
[{\mb{K}}_-,{\mb{K}}_+]{\mb{K}}^{p-1-q}_+=
2\sum_{q=0}^{p-1}{\mb{K}}^q_+{\mb{K}}_0{\mb{K}}^{p-1-q}_+,$$ 
which leads in the expression of $\mu_{n;m}$ to the term 
$2{\mb{K}}_0{\mb{K}}^{p-1}_+$,
and, after acting to the left with ${\mb{K}}_0$, we get  a 0 contribution.

Similarly, if $n>m$, then
$$[{\mb{K}}^p_-,{\mb{K}}_+]=
-\sum_{s=0}^{p-1}{\mb{K}}_-^s[{\mb{K}}_0,{\mb{K}}_-]{\mb{K}}_-^{p-s-1},$$
and $s=p-1$ in the sum. Acting on the right with ${\mb{K}}_0$, the
contribution is also 0 because of the action on the right with 
${\mb{K}}_-^{p-1}$.\hfill $\gata$

\subsection*{Acknowledgments}

 The author is thankful to
 the organizers of {\it 2$^{~\text{nd}}$ Operator Algebras
and Mathematical Physics Conference}, Sinaia, Romania, June 26 -- July
4, 2003, 
 and  of the 
 {\it XXIII Workshop  on Geometric methods in Physics},
   June 27 --  July 3 2004, Bia\l owie\.{z}a,  Poland for the opportunity to
report results on this subject and for the financial
support for attending the conferences. Discussions 
 with Professor S. Twareque Ali and
Professor Peter Kramer are kindly acknowledged. The author is grateful
to Professor Rolf Berndt and to  Dr. Adrian Tanas\u{a} for correspondence,
 to Professor  Karl-Hermann
Neeb for criticism and to Professor John Klauder  for his involved  interest.


\begin{thebibliography}{99}

\bibitem{ali}S. T. Ali, J.-P. Antoine and J.-P. Gazeau, Coherent states,
wavelets, and their generalizations, Springer-Verlag, New York (2000)

\bibitem{bar47}V. Bargmann, {\it Irreducible unitary representations of the
Lorentz group},  Ann. of Math.  {\bf 48} (1947) 568-640


\bibitem{bar}V. Bargmann, {\it On the Hilbert space of analytic functions
and the associated integral transform},  Commun. Pure
Appl. Math. {\bf 14}  (1961) 187-214
 
\bibitem{bar70}V. Bargmann, Group representations on Hilbert spaces of
analytic functions, in {\it Analytic methods in Mathematical Physics},
Edited by R.P Gilbert and R.G. Newton, Gordon and Breach, Science
publishers, New York, London, Paris (1970) 27-63

\bibitem{bate}H. Bateman, Higher transcendental functions, Volume 2,
Mc Graw-Hill book,  New  York (1958)

\bibitem{morse}S. Berceanu and   C. A. Gheorghe, {\it On the construction of
perfect Morse functions on compact manifolds of coherent states},
 J. Math. Phys. {\bf 28} (1987) 2899-2907

\bibitem{sbcag}S. Berceanu and A. Gheorghe,
  {\it On equations of motion on Hermitian
symmetric  spaces},
 J. Math. Phys.  {\bf 33 } (1992) 998-1007

\bibitem{sbl}S. Berceanu and L. Boutet de Monvel,
{\it  Linear dynamical systems, coherent
state manifolds, flows and matrix Riccati equation},  J. Math. Phys.
{\bf 34}  (1993) 2353-2371

\bibitem{sbm}S. Berceanu and M. Schlichenmaier, {\it Coherent state
embeddings,
 polar
divisors and Cauchy formulas},
   J. Geom. Phys. {\bf 34}  (2000) 336-358; arXiv:  math. DG/9903105

\bibitem{sbctim}  S. Berceanu and  A. Gheorghe,
{\it Linear Hamiltonians on homogeneous K\"ahler ma\-nifolds
 of coherent states}, An. Univ. Timi\c soara Ser. Mat.-Inform.
   Vol XXXIX (2001) 31-56,  Special Issue: Mathematics;
 arXiv:math.DG/0408254 


\bibitem{last}S. Berceanu and A. Gheorghe,
{\it  Differential operators on orbits of coherent states},
 Rom. Jour.
Phys. {\bf  48} (2003) 545-556; arXiv:  math.DG/0211054


\bibitem{cluj}S. Berceanu,  Geometrical  phases on hermitian
symmetric spaces,
in {\it  Recent Advances in Geometry and Topology},
Editted by  Dorin Andrica,
Paul A. Blaga,   Cluj University Press (2004)  83-98;
    arXiv: math.DG/0408233 

\bibitem{sin}S. Berceanu,   Realization of coherent state algebras
by differential operators,
in {\it Advances in Operator Algebras and Mathematical Physics},
Edited by  F. Boca, O. Bratteli, R. Longo, H. Siedentop,
The Theta Foundation, Bucharest (2005) 1-24;  
arXiv: math.DG/0504053

\bibitem{sbj}S. Berceanu, {\it A holomorphic representation of the
semidirect sum of symplectic and Heisenberg Lie algebras}, 
 J. Geom. Symmetry Phys. 5 (2006) 5-13;  S. Berceanu,   A holomorphic 
representation of Jacobi algebra in several 
dimensions,  in {\it Perspectives in Operator Algebra and Mathematical
  Physics}, Edited by  F. P. Boca, R. Purice, S. Stratila, 
  The Theta Foundation, Bucharest (2008)  1-25; arXiv math.DG/060404381



 
\bibitem{berezin}F. A. Berezin, {\it The general concept of quantization},
  Commun. Math. Phys. {\bf 40} (1975) 153-174


\bibitem{berezin1}F. A. Berezin, {\it Models of Gross-Neveu type are
  quantization of a Classical Mechanics with a nonlinear phase space},
 Commun. Math. Phys. {\bf 63} (1978)  131-153

\bibitem{bern}R.  Berndt,  Sur l'arithmétique du corps des
fonctions elliptiques de niveau $N$, in {\it Seminar on number theory,
Paris 1982--83},   Progr. Math., 51,
Birkh\"auser Boston, Boston, MA (1984) 21-32

\bibitem{bb}R. Berndt and S. B\"ocherer, {\it Jacobi forms and
discrete series representations of the Jacobi group}, Math. Z. {\bf 204} (1990)
13-44 


\bibitem{bs}R. Berndt and R.  Schmidt, Elements of the representation
theory of the Jacobi group, Progress in Mathematics, 163,  Birkh\"auser
Verlag, Basel (1998) 

\bibitem{berndt}R. Berndt, {\it Coadjoint orbits and representations
of the Jacobi group}, IHES/M/03/37,  preprint (2003)

\bibitem{bi}I. Bialynicki-Birula, {\it
Solutions of the equations of motion in classical and quantum
theories}, Ann. Phys. NY {\bf 67} (1971) 252-273

\bibitem{bied}L. C. Biedenharn, {\it The quantum group
${\rm{SU}}_q(2)$ and a $q$-analogue of the boson operators},
J. Phys. A: Math. Gen. {\bf 22} (1989) 4581-4588

\bibitem{bog}N. N.  Bogoliubov, {\it A contribution to the
 theory of super-fluidity}, Bull. Acad. Sci. USSR. [Izvestia
Akad. Nauk. SSSR] {\bf 11} (1947) 77-90

\bibitem{CGR1}M.~Cahen, S.~Gutt and J.~Rawnsley, {\it  Quantization
of K\"{a}hler manifolds 
  {I}: Geometric interpretation of Berezin's quantization},
J. Geom. Phys. {\bf 7}  (1990) 45-62

\bibitem{dix} J. Dixmier, {\it Sur les alg\'ebres de Weyl},
Bull. Soc. Math. France {\bf 96} (1968) 209-242 

\bibitem{dod} V. V. Dodonov,
 {\it `Nonclassical' states in quantum optics: a `squeezed' review of 
the first 75 years}, J. Opt. B.:  Quantum semiclassical
Opt. {\bf 4}, (2002) R1-R33


\bibitem{dr}P. D. Drummond and Z. Ficek, Editors, 
Quantum squeezing, Springer, Berlin (2004) 

\bibitem{ez} M. Eichler and D. Zagier, The theory of Jacobi forms,
Progress in Mathematics, 55,
Birkh\"auser, Boston, MA (1985)

\bibitem{fol}G. B. Folland,  Harmonic analysis in phase space, Princeton
University Press, Princeton, New Jersey (1989)

\bibitem{gl}R. J. Glauber, {\it Coherent and incoherent states of the radiation
field},  Phys. Rev. {\bf 131} (1963) 2766-2788
 
\bibitem{grad} I. S Grad\v ste\u\i n i  I. M. Ry\v zik, 
{\cyr Tablitsy integralov, summ, ryadov i proizvedeni\u\i }. 
(Russian) [Tables of integrals,
sums, series and products.] Gosudarstv. Izdat. Fiz.-Mat. Lit., 
Moscow (1963) 

\bibitem{hn}J. Hilgert, K-H. Neeb and  B. \O rsted, {\it Conal Heisenberg algebras
and associated Hilbert spaces}, J. reine angew Math. {\bf 474} (1976) 67-112

\bibitem{ho}J. N. Hollenhorst, {\it Quantum limits on resonant-mass
gravitational-wave detectors}, Phys. Rev. D {\bf 19} (1979) 1669-1679

\bibitem{hol}T. Holstein and H. Primakoff, {\it
Field dependence of the intrinsic domain magnetization of a
ferromagnet}, Phys. Rev. {\bf 58} (1940) 1098-1113 


 \bibitem{itzik}C. Itzykson, {\it Remarks on boson commutation rules},
Commun. Math. Phys. {\bf 4} (1967) 92-122  

\bibitem{jos}A. Joseph, {\it Commuting polynomials in quantum
canonical operators and realizations of Lie algebras},
J. Math. Phys. {\bf 13} (1972) 351-357

\bibitem{cal1} E.  K\"ahler, {\it  Die
Poincar\'e-Gruppe},
Rend. Sem. Mat. Fis. Milano 53 (1983) 359-390 

 \bibitem{cal2} E.  K\"ahler, {\it The Poincar\'e group}, in {\it Clifford
algebras and their applications in mathematical physics} (Canterbury,
1985),  NATO Adv. Sci. Inst. Ser. C Math. Phys. Sci., 183,
Reidel, Dordrecht, (1986) 265-272 

\bibitem{cal3} E.  K\"ahler,
{\it Raum-Zeit-Individuum},
Rend. Accad. Naz. Sci. XL Mem. Mat. (5)  16  (1992)
115-177 


\bibitem{cal} Erich K\"ahler: Mathematische Werke; Mathematical
Works, Edited by R. Berndt, O. Riemenschneider, 
  Walter de Gruyter, Berlin-New York (2003)

\bibitem{ken}E. H. Kennard, {\it Zur Quantenmechanik einfacher
Bewegungstypen}, Zeit. Phys. {\bf 44} (1927) 326-352

\bibitem{kir}A. Kirillov, \'El\'ements de la th\'eorie des
representations, Editions Mir, Moscou (1974)

\bibitem{kir1}A. A. Kirillov, Lectures on the orbit method, Graduate
studies in Mathematics 64, American Mathematical Society, Providence,
Rhode Island (2004)

\bibitem{kir2}A. A. Kirillov, {\it Merits and demerits of the orbit method},
Bull. Amer. Math. Soc. {\bf 36} (1999) 43-73

\bibitem{kli}A. U. Klimyk, {\it On position and momentum operators in
the  $q$-oscillator}, J. Phys. A {\bf 38} (2005) 4447-4458

\bibitem{koby}S. Kobayashi, {\it Irreducibility of certain unitary
representations},  J. Math. Soc. Japan {\bf 20} (1968) 638-642

\bibitem{lie} S. Lie, {\it Theorie der Transformationsgruppen},
Math. Ann. {\bf 16} (1880) 441-528

 \bibitem{lis}W. Lisiecki, {\it Coherent state representations. A survey},
  Rep. Math. Phys. {\bf 35} (1995) 327-358

\bibitem{lu}E. Y. C. Lu, {\it New coherent states of the
electromagnetic field}, Lett. Nuovo. Cimento {\bf 2} (1971) 1241-1244

\bibitem{mo}B. R. Mollow and R.J. Glauber, {\it Quantum theory of
parametric amplifications}:I, Phys. Rev. {\bf 160} (1967) 1076-1096
   
\bibitem{neeb94} K.-H. Neeb, {\it Realization of general unitary highest
weight  representations}, Preprint, Technische Hochschule Darmstadt {\bf
1662} (1994) 


\bibitem{neeb}K.-H. Neeb, Holomorphy and Convexity in Lie Theory, de
 Gruyter Expositions in Mathematics 28, Walter de Gruyter, Berlin-New
York (2000)

\bibitem{nieto}M. N. Nieto and D. R. Truax, {\it 
Holstein-Primakoff/Bogoliubov transformations and the multiboson
system},  Fortschr. Phys. {\bf 45} (1997) 145-156

\bibitem{PWJ} Z. Pasternak-Winiarski and J. Wojcieszynski, {\it Bergman
spaces and kernel for holomorphic vector bundles},  Demonstratio
mathematica {\bf XXX}, (1997)  199-214

\bibitem{perG}A. M.  Perelomov,  Generalized Coherent States and their
Applications, Springer, Berlin (1986)

\bibitem{ps} I. I. Pyatetskii-Shapiro,  Automorphic functions and
the geometry of classical domains, Gordon \& Breach, New
York-London-Paris (1969)

\bibitem{rom} M. Rausch de Traubenberg, M. J. Slupinski and A. Tanas\u{a},
{\it Finite-dimensional Lie subalgebras of the Weyl algebra},
arXiv:math.RT/0504224 v2

\bibitem{raw}J. H. Rawnsley, {\it Coherent states and K\"ahler manifolds},
  Quart. J. Math. Oxford   {\bf 28} (1977) 403-415

\bibitem{satake}I. Satake, Algebraic structures of symmetric domains,
Publ. Math. Soc. Japan {\bf 14},  Princeton Univ. Press (1980)

\bibitem{si}R. Simon, E.C.G. Sudarshan and N. Mukunda, {\it Gaussian
pure states in quantum mechanics and the symplectic group},
Phys. Rev. A. {\bf 37} (1988) 3028-3038

\bibitem{sz}A. Simoni and F. Zaccaria, {\it On realization of
semi-simple Lie algebras with quantum canonical variables}, Nuovo
Cimento A (10) {\bf 59} (1969) 280-292

\bibitem{siv}S. Sivakumar, {\it Studies on nonlinear coherent states},
J. Opt. B: Quantum Semiclass. Opt. {\bf 2} (2000) R61-R75
\bibitem{stol}P. Stoler P, {\it Equivalence classes of minimum uncertainty packets},
 Phys. Rev. D {\bf  1} (1970) 3217-3219; ---, II,  Phys. Rev. D
{\bf 4} (1971) 1925-1926
\bibitem{wa}D.F. Walls, {\it Squeezed states of light}, Nature {\bf
306}
(1983) 141-146
\bibitem{yu}H.P. Yuen, {\it Two-photon coherent states of the
radiation field}, Phys. Rev. A {\bf 13} (1976) 2226-2243
\bibitem{zh}W.-M Zhang, D. H. Feng and R. Gilmore, {\it Coherent
states: theory and some applications}, Rev. Mod. Phys. {\bf 62} (1990)
867-927
\end{thebibliography}
\end{document}